\documentclass[review,onefignum,onetabnum]{siamart171218}

\usepackage{lipsum}
\usepackage{mathrsfs}
\usepackage{amsfonts}
\usepackage{graphicx}
\usepackage{epstopdf}
\usepackage{algorithmic}
\usepackage{amssymb}
\usepackage{amsmath}

\ifpdf
  \DeclareGraphicsExtensions{.eps,.pdf,.png,.jpg}
\else
  \DeclareGraphicsExtensions{.eps}
\fi

\newsiamremark{hypothesis}{Hypothesis}
\crefname{hypothesis}{Hypothesis}{Hypotheses}
\newsiamthm{claim}{Claim}
\newtheorem{remark}{Remark}
\nolinenumbers

\headers{Parabolic cylinder and Weber equations}{T. M. Dunster}

\title{Uniform asymptotic expansions for solutions of the parabolic cylinder and Weber equations}

\author{T. M. Dunster\thanks{Department of Mathematics and Statistics, San Diego State University, 5500 Campanile Drive, San Diego, CA 92182-7720, USA. 
  (\email{mdunster@sdsu.edu}, \url{https://tmdunster.sdsu.edu}).}}

\usepackage{amsopn}

\makeatletter
\newcommand*{\addFileDependency}[1]{
  \typeout{(#1)}
  \@addtofilelist{#1}
  \IfFileExists{#1}{}{\typeout{No file #1.}}
}
\makeatother

\nolinenumbers
 
\ifpdf
\hypersetup{
  pdftitle={Uniform asymptotic expansions for solutions of the parabolic cylinder and Weber equations},
  pdfauthor={T. M. Dunster}
}

\begin{document}

\maketitle

\begin{abstract}
 Asymptotic expansions are derived for solutions of the parabolic cylinder and Weber differential equations. In addition the inhomogeneous versions of the equations are considered, for the case of polynomial forcing terms. The expansions involve exponential, Airy and Scorer functions and slowly varying analytic coefficient functions involving simple coefficients. The approximations are uniformly valid for large values of the parameter and unbounded real and complex values of the argument. Explicit and readily computable error bounds are either furnished or available for all approximations.
\end{abstract}



\section{Introduction} 
\label{sec1}

We study the homogeneous ($h(z)=0$) and inhomogeneous ($h(z)=z^{R}$, $R=0,1,2,\cdots$) versions of the parabolic cylinder equation
\begin{equation} \label{1.1}
\frac{d^{2}y}{dz^{2}}- \left(\frac{1}{4}z^{2} \pm a \right)y=h(z),
\end{equation}
as well as the closely related versions of the Weber's equation
\begin{equation} \label{1.2}
\frac{d^{2}y}{dz^{2}}+ \left(\frac{1}{4}z^{2} \pm a \right)y=h(z),
\end{equation}
where $a$ is positive and large, and $z$ lies in unbounded complex domains. By the superposition principal the inhomogeneous monomial term can be extended to an arbitrary polynomial. Both equations are characterized as having two turning points when $a$ is large, either both real or both imaginary. For a definition of turning points, see \cite[Sect. 2.8(i)]{NIST:DLMF}.

All solutions of the above equations can be expressed in terms of the parabolic cylinder function $U(a,z)$, which is a solution of homogeneous equation
\begin{equation}  \label{1.3}
\frac{d^{2}y}{dz^{2}}-\left(\frac{1}{4}z^{2}+a \right)y=0,
\end{equation}
and has the integral representations \cite[Eqs. 12.5.1 and 12.5.4]{NIST:DLMF}

\begin{equation} \label{1.4}
U(a,z)=\frac{e^{-\tfrac{1}{4}z^{2}}}{\Gamma
\left(\tfrac{1}{2}+a \right)}
\int_0^\infty {t^{a-\tfrac{1}{2}}
e^{-\tfrac{1}{2}t^{2}-zt}dt}
\quad \left(\Re(a)>-\tfrac{1}{2} \right),
\end{equation}
and
\begin{equation} \label{1.5}
U(-a,z)=\sqrt {\frac{2}{\pi }} e^{
\tfrac{1}{4}z^{2}}\int_0^\infty {t^{a-\tfrac{1}{2}}e^{-\tfrac{1}{2}t^{2}}\cos \left( {zt-\tfrac{1}{2}\pi a+\tfrac{1}{4}\pi } \right)dt} \quad \left(\Re(a)>-\tfrac{1}{2} \right).
\end{equation}
It can also be expressed in terms of the confluent hypergeometric function \cite[Sect. 13.2]{NIST:DLMF} by
\begin{equation} \label{1.6}
U(a,z)=2^{-\tfrac{1}{4}-\tfrac{1}{2}a}
e^{-\tfrac{1}{4}z^{2}} U\left( \tfrac{1}{2}a+\tfrac{1}{4},
\tfrac{1}{2},
\tfrac{1}{2}z^{2}
\right).
\end{equation}

The fundamental property of $U(a,z)$ is that it is the unique solution of (\ref{1.3}) that is recessive in the sector $|\arg(z)|\leq \pi /4$, whereas all other independent solutions are dominant in this sector. Specifically, as $z\to \infty$
\begin{equation} \label{1.7}
U(a,z)\sim z^{-a-\frac{1}{2}}
e^{-\tfrac{1}{4}z^{2}}
\quad \left( {\vert \arg (z)\vert \le \tfrac{3}{4}\pi -\delta } 
\right).
\end{equation}

Parabolic cylinder functions have a number of physical applications, most notably as solutions of the Helmholtz equation \cite[Sect. 12.17]{NIST:DLMF}. See also \cite{Buchholz:1969:CHF} for further properties in connection with physical applications. Solutions of the inhomogeneous Weber equation (\ref{1.2}) play a role in the study of flow through porous layers \cite{Zaytoon:2016:WID},  \cite{Nield:2009:TET}.

Mathematical applications of parabolic cylinder functions include the uniform asymptotic approximation of solutions of differential equations having two turning points, most significantly for the case where a pair can coalesce \cite{Olver:1975:SOL}. Parabolic cylinder functions are also used in approximating contour integrals having a coalescing saddle point and an algebraic singularity \cite[Chap. 22]{Temme:2015:AMF}. 

For numerical methods of evaluating parabolic cylinder and Weber functions for real $a$ and $z$ see \cite{Gil:2004:IRC}, \cite{Gil:2006:RPC}, \cite{Gil:2006:CRPC}, \cite{Gil:2011:A914}, \cite{Gil:2011:FAC}. In a subsequent paper it is intended to apply the new results in this paper to develop numerical algorithms to compute the parabolic cylinder functions having complex argument.

There have been many investigations into the asymptotic behaviour of parabolic cylinder functions, the most notable being by Olver \cite{Olver:1959:USE} who considered both real and complex values of $a$ and $z$. He used elementary functions in a Liouville Green (LG) expansion in intervals or domains free of turning points, and Airy functions in domains or intervals containing a turning point. Error bounds were not given and the coefficients are harder to compute near the turning points. Temme \cite{Temme:2000:NAA} modified the approximations found in \cite{Olver:1959:USE} by redefining certain integration constants so that the coefficients approach zero at infinity, and included useful numerical examples. In his paper only real variables were considered and there were no error bounds. Jones \cite{Jones:2006:PCF} did provide some error bounds, but these are quite limited and only real $z$ was considered.

In this paper we also obtain LG expansions involving elementary functions for the homogeneous equations away from turning points, and Airy expansions for the homogeneous equations in domains or intervals containing a turning point. In addition, we obtain expansions involving elementary functions for the inhomogeneous equations away from turning points, and expansions involving Scorer functions for solutions of the inhomogeneous equations in domains containing a turning point.

The significance of the present paper is that explicit and computable error bounds are available for all real and complex numerically satisfactory solutions of (\ref{1.1}) and (\ref{1.2}), and for the turning point cases our coefficients are easier to compute than those in the literature. Our new results come from the asymptotic theory provided in the recent papers \cite{Dunster:2020:ASI}, \cite{Dunster:2020:LGE} and \cite{Dunster:2020:SEB}. In \cite{Dunster:2020:LGE} LG expansions were derived, with error bounds, where the coefficients in the expansions appear in the exponent of an exponential function, rather than the typical case as a factor multiplying the exponential \cite[Chap. 10]{Olver:1997:ASF}. As a result the coefficients and error bounds of \cite{Dunster:2020:LGE} are easier to compute.

In \cite{Dunster:2020:SEB} simple error bounds were derived for Airy turning point expansions for homogeneous differential equations having a turning point in the complex plane. A feature is that the coefficients that appear in the asymptotic expansion of two slowly varying coefficient functions are very easy to compute, either directly if not too close to the turning point, or via Cauchy's integral theorem in a neighbourhood of the turning point. The same is true for the expansions furnished by \cite{Dunster:2020:ASI}, where new asymptotic solutions, with error bounds, were constructed for inhomogeneous differential equations having a turning point.

The plan of this paper is as follows. In \cref{sec2} we consider (\ref{1.1}) with the upper sign taken on $a$. In this case the two turning points are on the imaginary axis, and in the homogeneous case we only need to derive LG expansions, using the theory of \cite{Dunster:2020:LGE}. For the inhomogeneous case we obtain expansions that involve elementary functions using \cite{Dunster:2020:ASI}. 

In \cref{sec3} we consider (\ref{1.1}) with the lower sign taken on $a$. This time the turning points are on the real axis, and in the homogeneous case we obtain expansions involving Airy functions, using \cite{Dunster:2020:SEB}. For the inhomogeneous case we obtain expansions that involve Scorer functions, using \cite{Dunster:2020:ASI}.

In \cref{sec4,,sec5} we consider (\ref{1.2}) with the upper and lower signs taken on $a$, respectively. Similar expansions are derived as in the previous two sections. There are significant differences however, particularly in the behaviour of the solutions in different sectors, and the process of identification of standard solutions.

In all the inhomogeneous cases we define appropriate fundamental particular solutions that are not exponentially large at infinity in certain sectors of the complex plane. Using the new theory of \cite{Dunster:2020:ASI} uniform asymptotic approximations are constructed for these solutions, and connection formulas are provided to aid in the approximation of all other solutions.

Although the parameter $a$ is considered real here, our results can be extended to complex $a$ without significant complication, since all the asymptotic theories being applied are valid for large complex values of the parameter. The forms of the expansions would remain the same, but the regions of validity would be modified accordingly, as well as the error bounds. We do not pursue this since in most applications the parameter is real.

\section{Parabolic cylinder functions with positive parameter} 
\label{sec2}

Let $a=\tfrac{1}{2}u$, $z$ be replaced by $\sqrt{2u} z$, and $w(u,z)=y(\tfrac{1}{2}u,\sqrt {2u} z)$. Then (\ref{1.3}) becomes
\begin{equation} \label{2.1}
\frac{d^{2}w}{dz^{2}}-u^{2}(z^{2}+1)w=0.
\end{equation}

Our approximations, and their regions of validity, involve the following function
\begin{equation} \label{2.2}
\bar{\xi}= \int_{0}^{z}\left(t^2+1\right)^{1/2} dt
=\frac{1}{2}z\sqrt {z^{2}+1}+\frac{1}{2}\ln\left( z+\sqrt {z^{2}+1} \right),
\end{equation}
where and throughout bars do not denote complex conjugate, unless otherwise noted. The branch is chosen so that $\bar{\xi}$ is real when $z$ is, both being of the same sign, and by continuity elsewhere in the plane with cuts along $z=\pm iy$, $1\leq y < \infty$. In what follows certain paths of integration can pass across these cuts, with the understanding that, if so, $\bar{\xi}$ varies continuously as a function of $z$ on the path in question.

\begin{figure}[htbp]
  \centering
  \includegraphics[width=\textwidth,keepaspectratio]{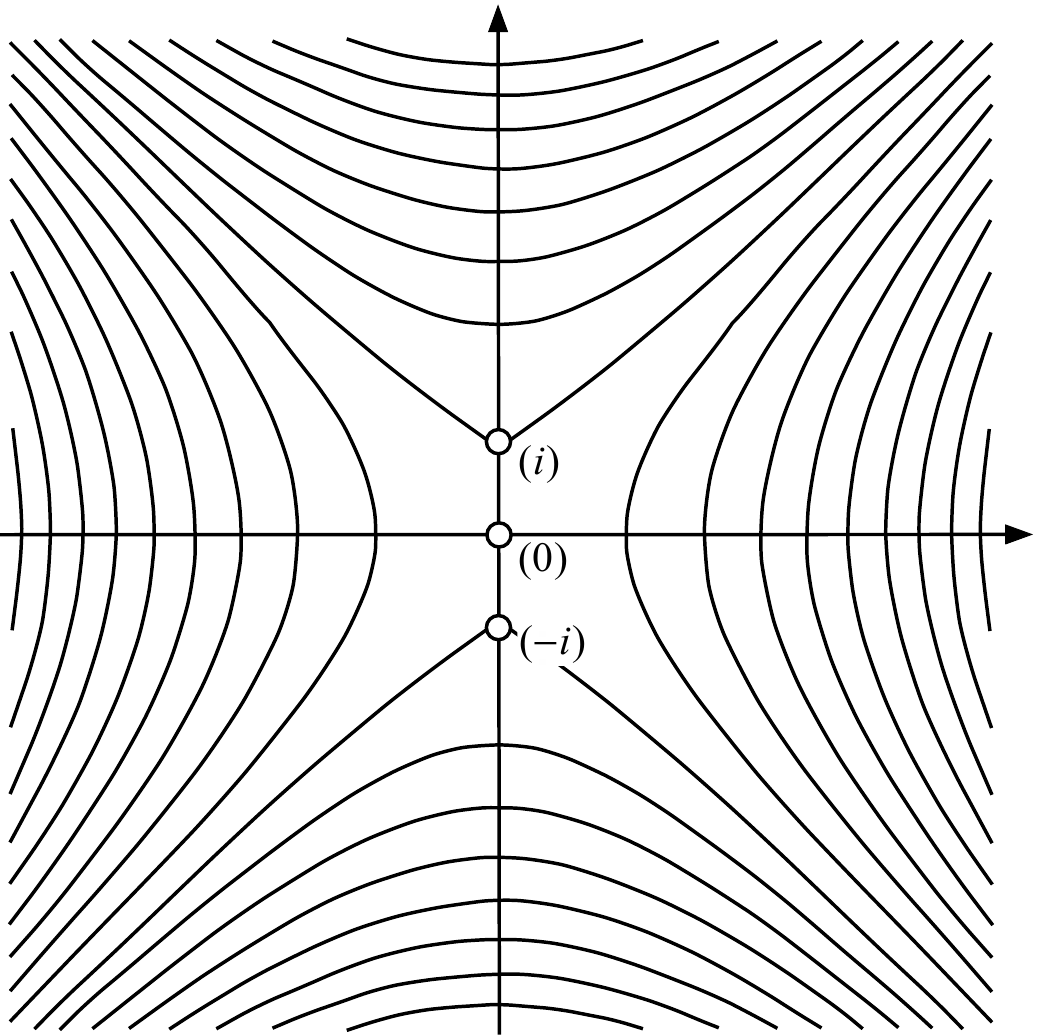}
  \caption{Level curves $\Re(\bar{\xi})=\mathrm{constant}$}
  \label{fig:fig1pcf}
\end{figure}

This variable comes from the Liouville transformation on (\ref{2.1}); see \cite[Chap. 10]{Olver:1997:ASF}. The regions of validity are based on the level curves $\Re(\bar{\xi})=\mathrm{constant}$, some of which are shown in \cref{fig:fig1pcf}. In this figure the points $z=\pm i$ are shown. These are turning points of (\ref{2.1}), and our LG approximations are not valid in their vicinity. However, this situation is covered in \cref{sec3}, where Airy function expansions are used to provide approximations which are valid at a turning point. Our main priority in this section is to obtain asymptotic approximations in unbounded complex domains that contain the whole real axis.

We also remark that the level curves emanating from $z=\pm i$ divide the plane into four regions, in each of which a solution of (\ref{2.1}) is recessive (and dominant elsewhere). These four fundamental solutions are given by $U(\tfrac{1}{2}u,\pm \sqrt {2u} z)$ and $U(-\tfrac{1}{2}u,\pm \sqrt {2u} iz)$. We shall obtain uniform asymptotic approximations for the former pair, as they are numerically satisfactory solutions in a domain containing the real axis.

Coefficients in our expansions are defined as follows. Firstly, identifying (\ref{2.1}) in the standard form $d^{2}w/dz^{2}=\{u^2f(z)+g(z)\}w(z)$ where $f(z)=z^2+1$ and $g(z)=0$ we have that the Schwarzian derivative \cite[Chap. 10, Eq. (2.03)]{Olver:1997:ASF}
\begin{equation} \label{2.3}
\Phi(z) =\frac{4f(z) {f}^{\prime \prime }(z) -5{f}^{\prime 2}(z)}{16f^{3}(z)}+\frac{g(z)}{f(z)},
\end{equation}
takes the form
\begin{equation} \label{2.4}
\Phi(z)=
\frac {2-3z^{2}}{ 4\left(z^{2}+1 \right) ^{3}}.
\end{equation}

Then following \cite[Eqs. (1.12) - (1.14)]{Dunster:2020:SEB} we define the set of coefficients 
\begin{equation} \label{2.5}
\hat{F}_{1}(z) ={\tfrac{1}{2}}\Phi(z) ,\ \hat{F}_{2}(z) =-{\tfrac{1}{4}}f^{-1/2}(z) {\Phi }^{\prime}(z),
\end{equation}
and 
\begin{equation} \label{2.6}
\hat{F}_{s+1}(z) =-\frac{1}{2}f^{-1/2}(z) \hat{{F}
}_{s}^{\prime }(z) -\frac{1}{2}\sum\limits_{j=1}^{s-1}{\hat{F}%
_{j}(z) \hat{F}_{s-j}(z) }\ \left( {s=2,3,4\cdots }
\right).
\end{equation}
The coefficients appearing in the following asymptotic expansions are then given by
\begin{equation} \label{2.7}
\hat{E}_{s}(z) =\int {\hat{F}_{s}(z)
f^{1/2}(z) dz}\ \left( {s=0,1,2,\cdots }\right).
\end{equation}

In order to facilitate these integrations we introduce a new variable by
\begin{equation} \label{2.8}
\bar{\beta}=\frac{z}{\sqrt{z^{2}+1}},
\end{equation}
from which
\begin{equation} \label{2.9}
\frac {d \bar{\beta}}{dz}= \left( 1-\bar{\beta}^{2} \right) ^{3/2}.
\end{equation}
We choose branch of the square root so that $\bar{\beta}$ is real and positive when $z$ is real and positive, and continuous in the $z$ plane having cuts at $z=\pm iy$ ($1\leq y < \infty$) on the imaginary axis. We note in passing that from (\ref{2.2}) and (\ref{2.8}) we have
\begin{equation} \label{2.13}
\bar{\xi}=\frac{\bar{\beta}}{2\left(1-\bar{\beta}^{2}\right)} +\frac{1}{4}
\ln \left(
\frac{\bar{\beta}+1}{1-\bar{\beta}}
\right).
\end{equation}

Now from (\ref{2.4}), (\ref{2.5}), (\ref{2.8}) and (\ref{2.9}) we find the first two coefficients are given by
\begin{equation} \label{2.10}
\bar{\mathrm{E}}_{1}(\bar{\beta})=\tfrac{1}{24}
\bar{\beta}
\left(6-5\bar{\beta}^{2}\right),
\end{equation}
\begin{equation} \label{2.11}
\bar{\mathrm{E}}_{2}(\bar{\beta})=
\tfrac{1}{16}\left(1-\bar{\beta}^{2}\right)^{2} 
\left(5\bar{\beta}^{2}-2\right),
\end{equation}
and from (\ref{2.7}) and (\ref{2.9}) for $s=2,3,4\cdots$
\begin{equation} \label{2.12}
\bar{\mathrm{E}}_{s+1}(\bar{\beta}) =
-\frac{1}{2} \left(1-\bar{\beta}^{2} \right)^{2}\bar{\mathrm{E}}_{s}^{\prime}(\bar{\beta})
-\frac{1}{2}\int_{\sigma(s)}^{\bar{\beta}}
\left(1-p^{2} \right)^{2}
\sum\limits_{j=1}^{s-1}
\bar{\mathrm{E}}_{j}^{\prime}(p)
\bar{\mathrm{E}}_{s-j}^{\prime}(p) dp,
\end{equation}
where $\sigma(s)$ is a suitably chosen constant. Our choice is $\sigma(s)=1$ for $s$ odd and $\sigma(s)=0$ for $s$ even. By induction it is easy to verify that each $\bar{\mathrm{E}}_{2s}(\bar{\beta})$ is even, and each $\bar{\mathrm{E}}_{2s+1}(\bar{\beta})$ is odd (and hence of course $\bar{\mathrm{E}}_{2s+1}(0)=0$).
Moreover, with our choice of $\sigma(s)$ we observe that $\bar{\mathrm{E}}_{2s}(1)=0$, which simplifies the expansions in subsequent sections.

We can now apply \cite[Thm. 1.1]{Dunster:2020:LGE} to obtain the solutions
\begin{equation} \label{2.14}
W_{1}(u,\bar{\beta}) =
\exp \left\{u\bar{\xi}
+\sum\limits_{s=1}^{n-1}{
\frac{
\bar{\mathrm{E}}_{s}(\bar{\beta}) 
-\bar{\mathrm{E}}_{s}(-1)}
{u^{s}}}
\right\} \left\{1+\bar{\eta}_{n,1}(u,\bar{\beta}) \right\},
\end{equation}
and 
\begin{equation} \label{2.15}
W_{2}(u,\bar{\beta}) =
\exp \left\{-u\bar{\xi}
+\sum\limits_{s=1}^{n-1}{(-1)^{s}
\frac{
\bar{\mathrm{E}}_{s}(\bar{\beta}) 
-\bar{\mathrm{E}}_{s}(1)}
{u^{s}}}
\right\} \left\{1+\bar{\eta}_{n,2}(u,\bar{\beta}) \right\}.
\end{equation}
Bounds on the error terms are given by
\begin{equation} \label{2.16}
\left\vert \bar{\eta}_{n,j}(u,\bar{\beta}) \right\vert \leq u^{-n}
\omega_{n,j}(u,\bar{\beta}) \exp \left\{u^{-1}\varpi_{n,j}(u,\bar{\beta}) +u
^{-n}\omega_{n,j}(u,\bar{\beta}) \right\}, 
\end{equation}
where 
\begin{multline} \label{2.17}
\omega_{n,j}(u,\bar{\beta}) =2\int_{(-1)^{j}}^{\bar{\beta}}{\left\vert
\bar{\mathrm{E}}_{n}^{\prime}(p) d p\right\vert } \\ 
+\sum\limits_{s=1}^{n-1}\frac{1}{u^{s}}
\int_{(-1)^{j}}^{\bar{\beta}}{\left\vert \sum\limits_{k=s}^{n-1}
(1-p^{2})^{2}
\bar{\mathrm{E}}_{k}^{\prime}(p)
\bar{\mathrm{E}}_{s+n-k-1}^{\prime}(p) dp
\right\vert },
\end{multline}
and 
\begin{equation} \label{2.18}
\varpi _{n,j}(u,\bar{\beta}) =4\sum\limits_{s=0}^{n-2}\frac{1}{
u^{s}}
\int_{(-1)^{j}}^{\bar{\beta}}
\left\vert
\bar{\mathrm{E}}_{s+1}^{\prime}(p)dp\right\vert .
\end{equation}

Here the paths of integration are taken along paths that avoid $z=\pm i$ (i.e. $\bar{\beta}$ unbounded) and on which $\Re(\bar{\xi})$ is monotonic. They must also consist of a finite chain of $R_{2}$ arcs (as defined in \cite[Chap. 5, sec. 3.3]{Olver:1997:ASF}).

\begin{figure}[htbp]
  \centering
  \includegraphics[width=\textwidth,keepaspectratio]{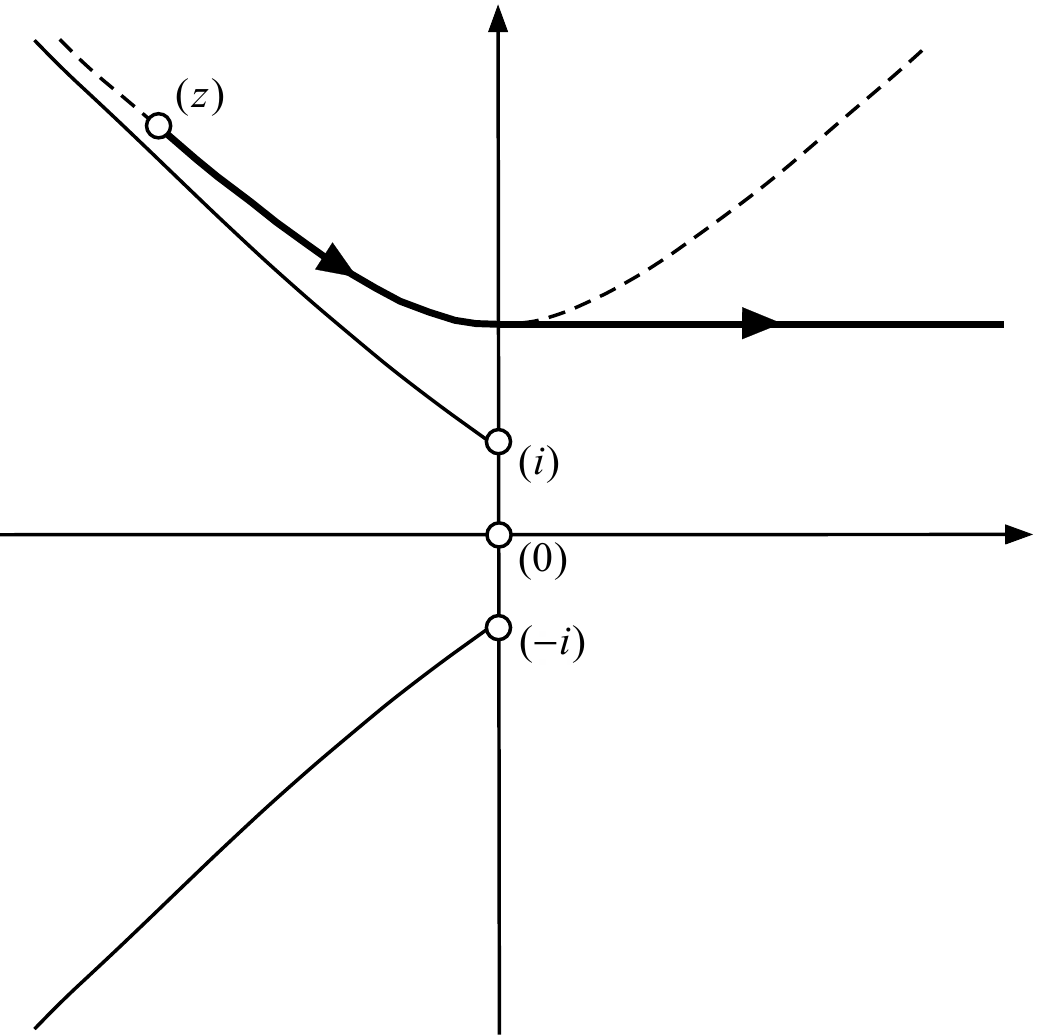}
  \caption{Path of integration}
  \label{fig:fig1aapcf}
\end{figure}

Consider $W_{2}(u,\bar{\beta})$, which is recessive as $\Re(z) \rightarrow \infty$, and the associated bounds (\ref{2.16}) - (\ref{2.18}) with $j=2$. An example of a suitable path in the $z$ plane is shown in \cref{fig:fig1aapcf} (the thick line) for a point lying in the second quadrant, and above the level curve emanating from $z=i$. The map of this path under (\ref{2.8}) is used for the integrals in (\ref{2.17}) and (\ref{2.18}).

In order for the error bound (\ref{2.16}) to be as sharp as possible such a path should be chosen to be as far away from the singularity $z=i$, but also maintaining the monotonicity requirement described above as $\bar{\xi}$ varies continuously along the path. Thus we choose it to consist of the union of a segment of the level curve passing through $z$, starting from this point and ending on the imaginary axis, along with a horizontal line extending to $\Re(z)=\infty$. The dashed line in the figure is the level curve in question. 

The region of validity for (\ref{2.15}) is the whole complex plane, with the exception of points lying on the two level curves emanating from $z = \pm i$ in the left half plane (see \cref{fig:fig1aapcf}). This is because all points except these can be linked to $\Re(z)=\infty$ by a path on which $\Re(\bar{\xi})$ is monotonic.

Note that as $z$ gets closer to the cut emanating from $z=i$ the path necessarily must get closer to the singularity $z=i$. This illustrates why all points on, and close, to this cut must be excluded, along with their conjugates. Similarly, the region of validity for (\ref{2.14}) is the reflection about the imaginary axis of the corresponding domain for (\ref{2.15}).

From the error bounds we note the doubly asymptotic behaviour that $\bar{\eta}_{n,2}(u,\bar{\beta})$ is $\mathcal{O}(u^{-n})$ as $u \rightarrow \infty$ in this unbounded domain, and also is $\mathcal{O}(z^{-2n})$ as $z \rightarrow \infty$ to the right of the above-mentioned level curves emanating from $z = \pm i$. This a feature for all the error terms in this paper.

Let us now match the asymptotic solutions with the parabolic cylinder function. From (\ref{2.2}) we find that as $\Re(z) \rightarrow \pm \infty$
\begin{equation} \label{2.19}
\bar{\xi}=
 \pm \tfrac{1}{2}z^{2} \pm \tfrac{1}{2}\ln(2z)
 \pm \tfrac{1}{4}
+\mathcal{O}(z^{-2}),
\end{equation}
and so from (\ref{1.7}) and (\ref{2.15}) and matching unique recessive solutions at $\Re(z)=\mp \infty$ we derive
\begin{equation} \label{2.21}
U\left(\tfrac{1}{2}u,-\sqrt {2u} z\right)
=\left( \frac{2e}{u} \right)^{u/4}
\frac{1}{\left\{2u
\left(1+z^{2}\right)
\right\}
^{1/4}}
W_{1}(u,\bar{\beta}),
\end{equation}
and
\begin{equation} \label{2.20}
U\left(\tfrac{1}{2}u,\sqrt {2u} z\right)
=\left( \frac{2e}{u} \right)^{u/4}
\frac{1}{\left\{2u
\left(1+z^{2}\right)
\right\}
^{1/4}}
W_{2}(u,\bar{\beta}).
\end{equation}

Although these can be differentiated, we find a simpler expansion with sharper bounds by looking at the the differential equation satisfied by the derivative. In particular, from (\ref{2.1}) we arrive at
\begin{equation} \label{2.24}
\frac{d^{2}Y}{dz^{2}}=\left\{
u^{2} (1+z^{2})
+\frac {2z^{2}-1}{(1+z^{2})^{2}}
\right\}Y,
\end{equation}
whose solutions are given by
\begin{equation} \label{2.25}
Y=(1+z^{2})^{-1/2}\partial w /\partial z.
\end{equation}

The analysis follows similarly to above (and using tildes instead of bars), so that in place of (\ref{2.4}), we have from (\ref{2.3})
\begin{equation} \label{2.26}
\tilde{\Phi}(z)=
\frac {5z^{2}-2}{ 4\left(z^{2}+1 \right) ^{3}},
\end{equation}
and hence similarly to (\ref{2.10}) - (\ref{2.12}) the coefficients
\begin{equation} \label{2.27}
\tilde{\mathrm{E}}_{1}(\bar{\beta})=\tfrac{1}{24}
\bar{\beta}
\left(7\bar{\beta}^{2}-6\right),
\end{equation}
and
\begin{equation} \label{2.28}
\tilde{\mathrm{E}}_{2}(\bar{\beta})=
\tfrac{1}{16}\left(1-\bar{\beta}^{2}\right)^{2} 
\left(2-7\bar{\beta}^{2}\right),
\end{equation}
with subsequent coefficients also given by (\ref{2.12}) (with bars replaced by tildes). Then on matching in a similar manner to (\ref{2.20}) and (\ref{2.21}) we arrive at our desired expansions
\begin{multline} \label{2.30}
U'\left(\tfrac{1}{2}u,-\sqrt {2u} z\right)=
-\frac{1}{2}
\left(\frac{2e}{u} \right)^{u/4}
\left\{2u\left(1+z^{2}\right)
\right\}^{1/4} \\
\times
\exp \left\{u\bar{\xi}
+\sum\limits_{s=1}^{n-1}
\frac{
\tilde{\mathrm{E}}_{s}(\bar{\beta}) 
-\tilde{\mathrm{E}}_{s}(-1)}{u^{s}}
\right\}
\left\{1+\tilde{\eta}_{n,1}(u,z) \right\},
\end{multline}
and
\begin{multline} \label{2.29}
U'\left(\tfrac{1}{2}u,\sqrt {2u} z\right)=
-\frac{1}{2}
\left(\frac{2e}{u} \right)^{u/4}
\left\{2u\left(1+z^{2}\right)
\right\}^{1/4}
 \\
\times
\exp \left\{-u\bar{\xi}
+\sum\limits_{s=1}^{n-1}(-1)^{s}
\frac{\tilde{\mathrm{E}}_{s}(\bar{\beta}) 
-\tilde{\mathrm{E}}_{s}(1)}{u^{s}}
\right\}
\left\{1+\tilde{\eta}_{n,2}(u,z) \right\},
\end{multline}
where $U'(a,z)=dU(a,z)/dz$. These are valid in the same domains as for (\ref{2.21}) and (\ref{2.20}) respectively, and the error terms have the same bounds (\ref{2.16}) - (\ref{2.18}), but with $\bar{\mathrm{E}}$ replaced by $\tilde{\mathrm{E}}$.

\subsection{Inhomogeneous equation}
Let us now consider the equation
\begin{equation} \label{2.32}
\frac{d^{2}y}{dz^{2}}-\left(\frac{1}{4}z^{2}+a \right)y=z^{R} \quad (R=0,1,2,\cdots).
\end{equation}
For $j \in \{0,1,2,3\}$, $k \in \{1,2,3\}$ with $j<k$ we have by variation of parameters the particular solutions
\begin{multline} \label{2.33}
 U_{R}^{(j,k)}(a,z)
 =\frac{1}
{\mathscr{W}\{U_{j}(a,z),U_{k}(a,z)\}}
\left[
U_{k}(a,z)
 \int_{i^{j}\infty}^z \right.
 t^{R}U_{j}(a,t)dt \\
\left. -U_{j}(a,z)
 \int_{i^{k}\infty}^z 
 t^{R}U_{k}(a,t)dt
 \right],
\end{multline}
where for $j=0,1,2,3$
\begin{equation} \label{2.34}
U_{j}(a,z)
=U\left((-1)^{j}a,(-i)^{j}z\right).
\end{equation}

Using the Wronskian \cite[Eq. 12.2.11]{NIST:DLMF}
\begin{equation} \label{2.35}
\mathscr{W}\{U(a,z),U(a,-z)\}
=\frac{\sqrt{2\pi}}
{\Gamma\left(a+\tfrac{1}{2}\right)},
\end{equation}
we have the particular solution we are most interested in
\begin{equation} \label{2.36}
 U_{R}^{(0,2)}(a,z)
 =-\frac{\Gamma\left(a+\tfrac{1}{2}\right)}
{\sqrt{2\pi}}\left[U(a,z)
 \int_{-\infty}^z t^{R}U(a,-t)dt
+U(a,-z) \int_{z}^{\infty}
 t^{R}U(a,t))dt \right].
\end{equation}

This is the unique solution which is not exponentially large as $z \rightarrow \pm \infty$, and in a domain which includes the whole real axis. It is straightforward to verify that $U_{R}^{(0,2)}(a,-z)=(-1)^{R}U_{R}^{(0,2)}(a,z)$.

We shall also consider connection formulas, and some of these involve $U_{R}^{(0,1)}(a,z)$, which from using \cite[Eq. 12.2.12]{NIST:DLMF} 
\begin{equation} \label{2.36a}
\mathscr{W}\left\{U(a,z),U(-a,-iz)\right\}
=e^{- (\frac{1}{2}a-\frac{1}{4})\pi i},
\end{equation}
is given by
\begin{multline} \label{2.36b}
 U_{R}^{(0,1)}(a,z)
 =e^{i\pi(\frac{1}{2}a-\frac{1}{4})}\left[
U(-a,-iz) \int_{\infty}^z \right.
 t^{R}U(a,t)dt \\
\left. -U(a,z) \int_{i\infty}^z 
 t^{R}U(-a,-it)dt
 \right].
\end{multline}

To obtain asymptotic expansions let
\begin{equation} \label{2.37}
w(u,z)=
(2u)^{-\frac{1}{2}R-1}y\left(\tfrac{1}{2}u,\sqrt{2u} z\right),
\end{equation}
where again $u=2a$. Then (\ref{2.32}) becomes
\begin{equation} \label{2.38}
\frac{d^{2}w}{dz^{2}}-u^{2}(z^{2}+1)w=z^{R}.
\end{equation}
We now can apply \cite[Thm. 4]{Dunster:2020:ASI} to obtain the desired expansions
\begin{equation} \label{2.39}
U_{R}^{(j,k)}\left(\tfrac{1}{2}u,\sqrt{2u} z\right)=
(2u)^{\frac{1}{2}R+1}\bar{w}^{(j,k)}_{R}(u,z),
\end{equation}
where
\begin{equation} \label{2.40}
\bar{w}^{(j,k)}_{R}(u,z)=\frac{1}{u^{2}}\sum\limits_{s=0}^{n-1} \frac{\bar{G}_{s,R}(z)}{u^{2s}} +\bar{\varepsilon}_{n,R}^{(j,k)}(u,z),
\end{equation}
in which
\begin{equation} \label{2.41}
\bar{G}_{0,R}(z)=-z^{R}/\left(z^{2}+1\right),
\end{equation}
and
\begin{equation} \label{2.42}
\bar{G}_{s+1,R}(z)=\bar{G}''_{s,R}(z)/\left(z^{2}+1\right)
\quad (s=0,1,2,\cdots ).
\end{equation}

The error terms $\bar{\varepsilon}_{n,R}^{(j,k)}(u,z)$ are $\mathcal{O}(u^{-2n-2})$ uniformly in certain unbounded domains, as described below. In these domains they satisfy the bounds
\begin{multline} \label{2.43}
 \left|\bar{\varepsilon}_{n,R}^{(j,k)}(u,z) \right|
 \le \frac{1}{u^{2n+2}}  \Bigg\{
 \left| \bar{G}_{n,R}(z)\right| \Bigg. \\
\left. +\frac{1}{2 |z^{2}+1|^{1/4}}
\int_{\bar{\mathcal{L}}^{(j,k)}(z)} {
\left| \left\{\left(t^{2}+1\right)^{1/4} \bar{G}_{n,R}(t)\right\}'dt
\right|}  \right\} \\ 
 +\frac{\bar{{L}}_{n,R}^{(j,k)}(z)}{8 u^{2n+3}
 {|z^{2}+1|}^{1/4} }\left\{ 
1-\frac{1}{8u }
\int_{\bar{\mathcal{L}}^{(j,k)}(z)} {\left| 
\frac {2-3t^{2}}{ \left(t^{2}+1 \right) ^{5/2}}dt 
\right|}  \right\}^{-1}
\\ \times
\int_{\bar{\mathcal{L}}^{(j,k)}(z)} {\left| 
\frac {2-3t^{2}}{ \left(t^{2}+1 \right) ^{5/2}} 
dt \right|},
\end{multline}
in which
\begin{multline} \label{2.44}
\bar{{L}}_{n,R}^{(j,k)}(z)=\sup _{t\in \bar{\mathcal{L}}^{(j,k)}(z)}
\left|\left(t^{2}+1\right)^{1/4} \bar{G}_{n,R}(t)  \right|
\\
+\frac{1}{2}\int_{\bar{\mathcal{L}}^{(j,k)}(z)} \left| 
\left\{\left(t^{2}+1\right)^{1/4}
\bar{G}_{n,R}(t)\right\}'dt
\right|.
\end{multline}
The path of integration $\bar{\mathcal{L}}^{(j,k)}(z)$ is an $R_{2}$ path that runs from $z=i^{j}\infty$ to $z=i^{k}\infty$ and passes through $z$, and avoids the singularities $z=\pm i$. As $t$ runs on these paths from one end point to the other $\Re(\bar{\xi})$ varies continuously and is monotonic. Since $\bar{G}_{n,R}(z) = \mathcal{O}(z^{R-2-4n})$ as $z \rightarrow \infty$ we assume $n>\frac{1}{4}R-\frac{3}{8}$ to ensure convergence of the integrals in (\ref{2.43}) and (\ref{2.44}).

The collection of all such points $z$ for which such a path exists defines the regions of validity. We denote these domains by $Z^{(j,k)}$. In general such domains depend on $\arg(u)$, but since we are only considering positive $u$ there is no such dependence in this paper.

\begin{figure}[htbp]
  \centering
  \includegraphics[width=\textwidth,keepaspectratio]{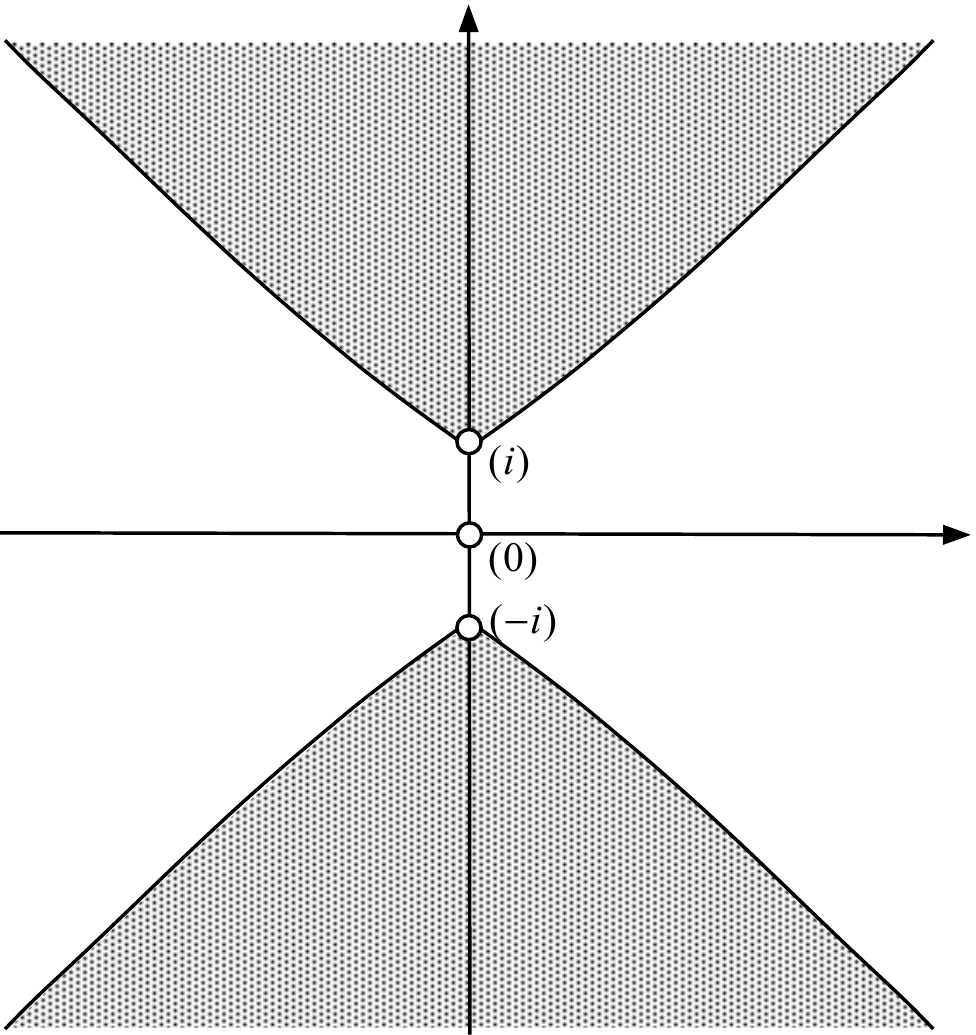}
  \caption{Domain $Z^{(0,2)}$}
  \label{fig:fig1apcf}
\end{figure}

\begin{figure}[htbp]
  \centering
  \includegraphics[width=\textwidth,keepaspectratio]{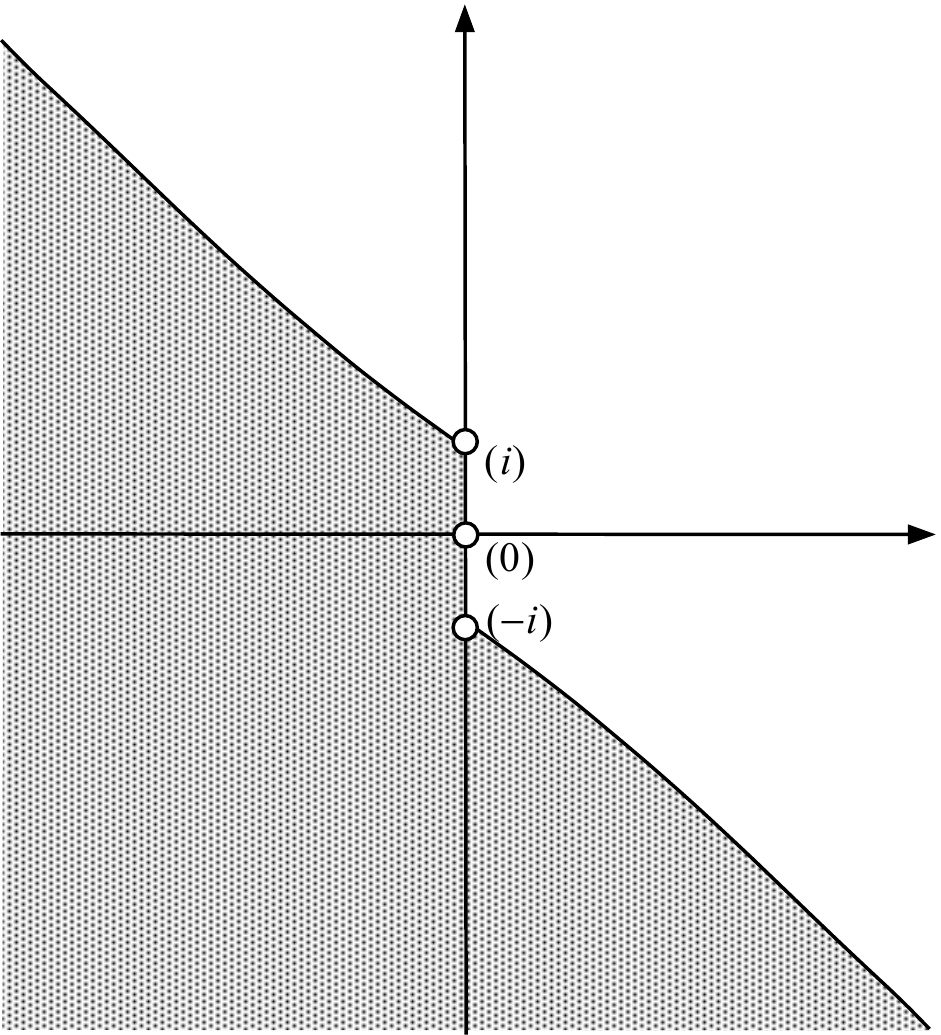}
  \caption{Domain $Z^{(0,1)}$}
  \label{fig:fig1bpcf}
\end{figure}

The domains $Z^{(0,2)}$ and $Z^{(0,3)}$ are depicted in \cref{fig:fig1apcf,,fig:fig1bpcf}, respectively: in these and subsequent figures shaded regions are excluded. Note that $Z^{(0,2)}$ contains the real $z$ axis, in accord with our previous comments about $U_{R}^{(0,2)}(a,z)$. The boundaries of these two domains are level curves emanating from the singularities, as well as a finite part of the imaginary axis in the case of $Z^{(0,1)}$. On inspecting the level curves in \cref{fig:fig1pcf} we deduce that all points on these boundaries must be excluded from the domains, i.e. they are open.

For the other cases we note that $Z^{(0,3)}$ is the complex conjugate of $Z^{(0,1)}$, and $Z^{(1,2)}$ (respectively $Z^{(2,3)}$) is the reflection of $Z^{(0,1)}$ (respectively $Z^{(0,3)}$) about the imaginary axis. 

We also remark that there is no appropriate path linking $z=\pm i \infty$, and hence $Z^{(1,3)}$ is empty. Thus (\ref{2.39}) and (\ref{2.40}), along with the accompanying error bounds, cannot be used for $j=1$, $k=3$. Instead appropriate connection formulas (given below) between the particular solutions must be used.

Let us now consider the aforementioned connection formulas between the above particular solutions. The relations we derive next will be useful in subsequent sections, as well as in their own right.

\begin{theorem}
\begin{equation} \label{2.47}
U_{R}^{(0,2)}(a,z)
=U_{R}^{(0,1)}(a,z)
+\Lambda_{R}(a)U(a,z),
\end{equation}
where for $a+\frac{1}{2}\neq 0,-1,-2,\cdots$
\begin{multline} \label{2.55}
\Lambda_{R}(a)=
2^{-\frac{1}{2}a+\frac{3}{2}R+\frac{1}{4}} 
\sqrt {\pi } e^{\left(\frac{1}{2}a
+\frac{1}{2}R-\frac{3}{4}\right)\pi i}
\Gamma \left(a+\tfrac{1}{2} \right) \\
\times
\mathbf{F}\left(\tfrac{1}{2}-\tfrac{1}{2}R,
-\tfrac{1}{2}R;
\tfrac{1}{2}a-\tfrac{1}{2}R+\tfrac{3}{4};\tfrac{1}{2}
\right),
\end{multline}
in which $\mathbf{F}$ is Olver's scaled hypergeometric function \cite[Eq. 15.2.2]{NIST:DLMF}
\begin{equation} \label{2.46}
\mathbf{F}(a,b;c,z)
=\frac{F(a,b;c,z)}{\Gamma(c)}
=\sum_{s=0}^{\infty}
\frac{(a)_{s}(b)_{s}z^{s}}{\Gamma(c+s)s!}.
\end{equation}
Moreover
\begin{multline} \label{2.56}
\Lambda_{R}(-a)=
-\left\{\tan(\pi a)+(-1)^{R}\sec(\pi a)+i\right\}
2^{-\frac{1}{2}a+\frac{3}{2}R-\frac{1}{4}} 
\pi \\
\times
\mathbf{F}\left(\tfrac{1}{2}-\tfrac{1}{2}R,
-\tfrac{1}{2}R;
\tfrac{1}{2}a-\tfrac{1}{2}R+\tfrac{3}{4};\tfrac{1}{2}
\right),
\end{multline}
provided $a \neq N+\tfrac{1}{2}$ for integer $N$ with $N-R$ even.
\end{theorem}

\begin{remark}
The hypergeometric function appearing in (\ref{2.55}) and (\ref{2.56}) consists of a finite number of terms. This is because either $\frac{1}{2}-\frac{1}{2}R$ or $-\frac{1}{2}R$ is a negative integer, and hence from \cite[Eq. 15.2.4]{NIST:DLMF} the series terminates. Also, for large $a$ the terms decrease rapidly, so the sum of the first few terms may suffice in obtaining a good approximation.
\end{remark}

\begin{proof}

We begin with the identity \cite[Eq. 2.11.2.1]{Prudnikov:1986:IAS}
\begin{multline} \label{2.45}
\int_{0}^{\infty}t^{R}U(a,t)dt=h_{R}(a)
:=2^{-\tfrac{1}{2}a-\tfrac{1}{2}R-\tfrac{3}{4}}
\sqrt{\pi}R!
\\ \times
\mathbf{F}\left(
\tfrac{1}{2}R+\tfrac{1}{2}, \tfrac{1}{2}R+1;
\tfrac{1}{2}a+\tfrac{1}{2}R+\tfrac{5}{4};\tfrac{1}{2}
\right).
\end{multline}

Now we know for some constant $\Lambda_{R}(a)$ that $U_{R}^{(0,1)}(a,z)$ and $U_{R}^{(0,2)}(a,z)$ are related by (\ref{2.47}). The reason is that both are particular solutions of (\ref{2.32}) that are not exponentially large at $z=+\infty$, and therefore their difference must be a multiple of $U(a,z)$, which is the unique solution of the associated homogeneous equation (\ref{1.3}) which is also not exponentially large at $z=+\infty$.

To find the constant we plug in $z=0$ in (\ref{2.36}), and refer to (\ref{2.45}), to obtain
\begin{equation} \label{2.48}
 U_{R}^{(0,2)}(a,0)
 =-\left\{1+(-1)^{R}\right\}
(2\pi)^{-1/2}\Gamma\left(a+\tfrac{1}{2}\right)
U(a,0)h_{R}(a).
\end{equation}
Next, we have from (\ref{2.33}), (\ref{2.34}), (\ref{2.36a}) and (\ref{2.45})
\begin{equation} \label{2.50}
 U_{R}^{(0,1)}(a,0)
 =e^{\left(\frac{1}{2}a-\frac{1}{4}\right)\pi i}
\left[i^{R+1}U(a,0)h_{R}(-a)
-U(-a,0)h_{R}(a)
\right].
\end{equation}
Now from \cite[12.2.17]{NIST:DLMF} we have
\begin{equation} \label{2.51}
U(-a,0)=(2/\pi)^{1/2}
\cos\left\{\left(\tfrac{1}{2}a-\tfrac{1}{4}\right)\pi
\right\}
\Gamma\left(a+\tfrac{1}{2}\right)U(a,0).
\end{equation}
And so from (\ref{2.47}), (\ref{2.48}) and (\ref{2.51})
\begin{multline} \label{2.52}
\Lambda_{R}(a)=
\left\{e^{\left(a-\frac{1}{2}\right)\pi i}-(-1)^{R}\right\}
 (2\pi)^{-1/2}\Gamma\left(a+\tfrac{1}{2}\right)
h_{R}(a)
\\
-e^{\left(\frac{1}{2}a
+\frac{1}{2}R+\frac{1}{4}\right)\pi i}h_{R}(-a).
\end{multline}

The function $h_{R}(-a)$ is of course given by (\ref{2.45}) with $a$ replaced by $-a$. We can simplify this by utilising the identity \cite[Eqs. 5.5.3 and 15.10.17]{NIST:DLMF}
\begin{multline} \label{2.53}
\frac{\mathbf{F}\left(a,b;-c;\tfrac{1}{2}\right)}
{\Gamma(a+c+1)\Gamma(b+c+1)}
=\frac{2^{a+b-1}\mathbf{F}\left(1-a,1-b;c+2;
\tfrac{1}{2}\right)}
{\Gamma(a)\Gamma(b)} \\
-\frac{1}{\pi}\sin(\pi c)
\mathbf{F}\left(a,b;a+b+c+1;\tfrac{1}{2}\right).
\end{multline}
Therefore, using this along with (\ref{2.45}) and \cite[Eq. 5.5.5]{NIST:DLMF} we obtain
\begin{multline} \label{2.54}
\frac{h_{R}(-a)}{\Gamma \left(a+\tfrac{1}{2} \right)} =
\sqrt{\frac{2}{\pi}}
\sin\left\{\left(\tfrac{1}{2}a
-\tfrac{1}{2}R-\tfrac{1}{4}\right)\pi\right\}h_{R}(a)
\\
+2^{-\tfrac{1}{2}a+\tfrac{3}{2}R+\tfrac{1}{4}} 
\sqrt {\pi } 
\mathbf{F}\left(\tfrac{1}{2}-\tfrac{1}{2}R,
-\tfrac{1}{2}R;
\tfrac{1}{2}a-\tfrac{1}{2}R+\tfrac{3}{4};\tfrac{1}{2}
\right).
\end{multline}
Then from (\ref{2.45}), (\ref{2.52}) and (\ref{2.54}) we arrive at (\ref{2.55}) and (\ref{2.56}).
\end{proof}

\section{Parabolic cylinder functions with negative parameter} 
\label{sec3}

Here we consider
\begin{equation} \label{3.1}
\frac{d^{2}y}{dz^{2}}-\left(\frac{1}{4}z^{2}-a \right)y=0,
\end{equation}
where again $a$ is large and positive. The fundamental solution is clearly $y(a,z)=U(-a,z)$. As before we let $a=\tfrac{1}{2}u$, have $z$ replaced by $\sqrt{2u} z$, and define $w(u,z)=y(\frac{1}{2}u,\sqrt {2u} z)$. Then (\ref{3.1}) becomes
\begin{equation} \label{3.3}
\frac{d^{2}w}{dz^{2}}-u^{2}(z^{2}-1)w=0.
\end{equation}

\begin{figure}[htbp]
  \centering
  \includegraphics[width=\textwidth,keepaspectratio]{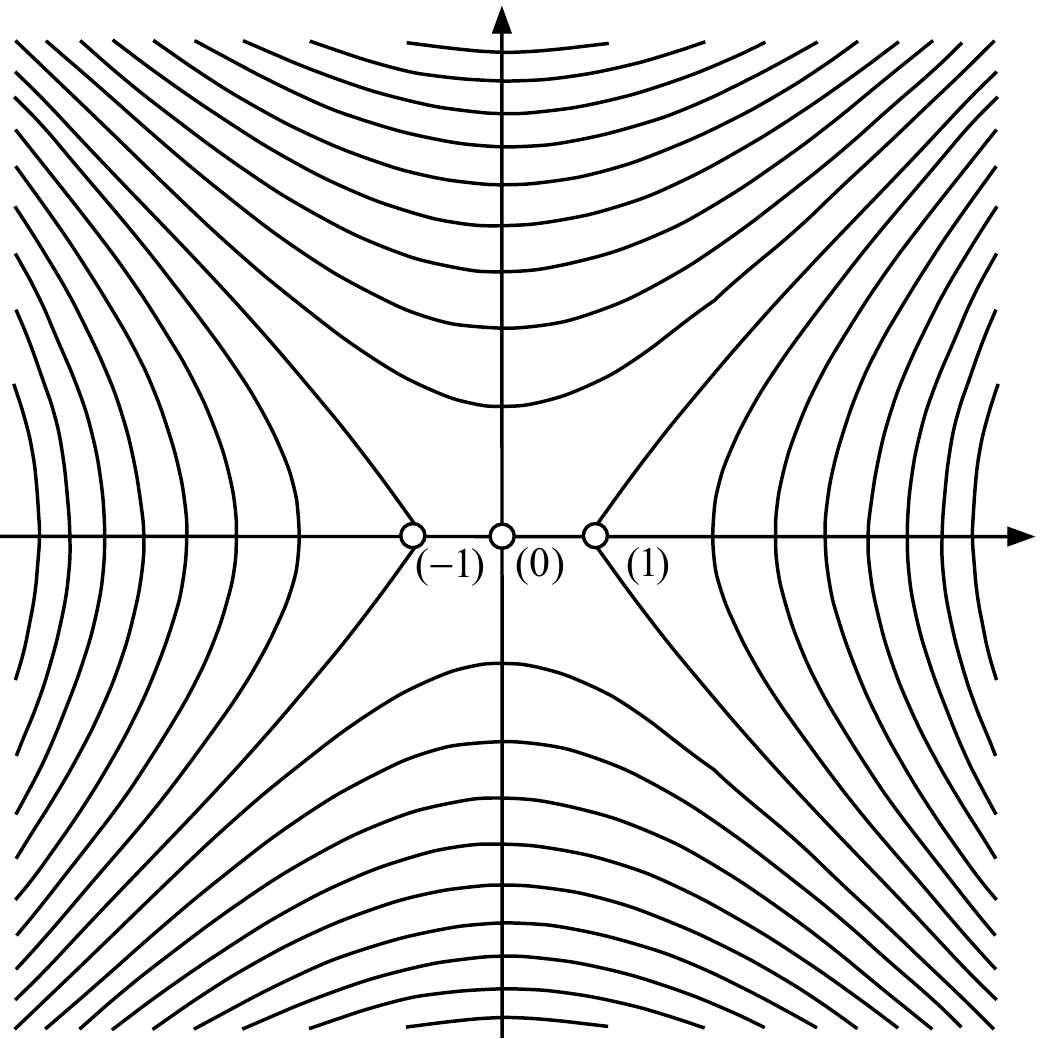}
  \caption{Level curves $\Re(\xi)=\mathrm{constant}$}
  \label{fig:fig2pcf}
\end{figure}

In place of (\ref{2.2}) we have instead the variables associated with LG ($\xi$) and Airy expansions ($\zeta$) given by
\begin{equation} \label{3.4}
\xi =\frac{2}{3}\zeta^{3/2}= \int_{1}^{z}\left(t^2-1\right)^{1/2} dt=\frac{1}{2}z\sqrt {z^{2}-1}-\frac{1}{2}\ln\left( z+\sqrt {z^{2}-1} \right).
\end{equation}
For $\xi$ we introduce a branch cut along the interval $(-\infty,1]$, and for $\zeta$ we introduce a branch cut along the interval $(-\infty,-1]$. With these cuts the branches are chosen so that $\xi \geq 0$ and $\zeta \geq0$ for $1\leq z<\infty$, and by continuity elsewhere in the respective cut $z$ plane. Some level curves in the $z$ plane corresponding to $\Re(\xi)=\mathrm{constant}$ are depicted in \cref{fig:fig2pcf}.

The turning point $z=1$ of (\ref{3.3}) is mapped to $\zeta =\xi=0$, and $\zeta$ is an analytic function of $z$ at this point. For $-1<z<1$ we note that $\zeta$ is negative and that
\begin{equation} \label{3.5}
\tfrac{2}{3}(-\zeta)^{3/2}
=\tfrac{1}{2}\arccos(z)-\tfrac{1}{2}z\sqrt{1-z^{2}}.
\end{equation}

The construction of the coefficients in our expansions follows similarly to \cref{sec2}, and in place of (\ref{2.3}) we have
\begin{equation} \label{3.6}
\Phi(z)=
-\frac {3z^{2}+2}{ 4\left(z^{2}-1 \right) ^{3}}.
\end{equation}
We let
\begin{equation} \label{3.7}
\beta=\frac{z}{\sqrt{z^{2}-1}},
\end{equation}
where the branch of the square root is positive for $z>1$ and is continuous in the plane having a cut along $[-1,1]$. Thus $\beta \rightarrow 1$ as $z \rightarrow \infty$ in any direction.

Similarly to (\ref{2.10}) - (\ref{2.12}) obtain the LG coefficients for this case as
\begin{equation} \label{3.8}
\mathrm{E}_{1}(\beta)=\tfrac{1}{24}\beta
\left(5\beta^{2}-6\right),
\end{equation}
\begin{equation} \label{3.9}
\mathrm{E}_{2}(\beta)=
\tfrac{1}{16}\left(\beta^{2}-1\right)^{2} 
\left(5\beta^{2}-2\right),
\end{equation}
and for $s=2,3,4\cdots$
\begin{equation} \label{3.10}
\mathrm{E}_{s+1}(\beta) =
\frac{1}{2} \left(\beta^{2}-1 \right)^{2}\mathrm{E}_{s}^{\prime}(\beta)
+\frac{1}{2}\int_{\sigma(s)}^{\beta}
\left(p^{2}-1 \right)^{2}
\sum\limits_{j=1}^{s-1}
\mathrm{E}_{j}^{\prime}(p)
\mathrm{E}_{s-j}^{\prime}(p) dp,
\end{equation}
where again $\sigma(s)=1$ for $s$ odd and $\sigma(s)=0$ for $s$ even, so that the even and odd coefficients are even and odd functions of $\beta$, respectively, with $\mathrm{E}_{2s}(1)=0$ ($s=1,2,3,\cdots$). 

The choice of $\sigma(s)$ for the even coefficients is for convenience in our identification of asymptotic solutions with parabolic cylinder functions. However, this choice of the lower integration limit for the odd coefficients is required for our subsequent Airy expansions to be valid in a full neighbourhood of the turning point; see \cite[Thm. 2.1]{Dunster:2017:COA}. From this reference we find that the integration constants for the odd coefficients must be chosen so that $(1-z)^{1/2}\mathrm{E}_{2s+1}(\beta) $ is a meromorphic function of $z$ at the turning point $z=1$. From (\ref{3.7}) we see that this is equivalent to $\mathrm{E}_{2s+1}(\beta)$ ($s=0,1,2,\cdots$) needing to be odd functions of $\beta$.

It is straightforward to verify that 
\begin{equation} \label{3.13}
\mathrm{E}_{s}(\beta)=(-1)^{s}
\bar{\mathrm{E}}_{s}(\beta),
\end{equation}
where $\bar{\mathrm{E}}_{s}$ are the coefficients used in the previous section. We also note from (\ref{3.4}) and (\ref{3.7}) that
\begin{equation} \label{3.14}
\xi=\frac{\beta}{2\left(\beta^{2}-1\right)}
 +\frac{1}{4}\ln \left(
\frac{\beta-1}{\beta+1}
\right).
\end{equation}

\begin{figure}[htbp]
  \centering
  \includegraphics[width=\textwidth,keepaspectratio]{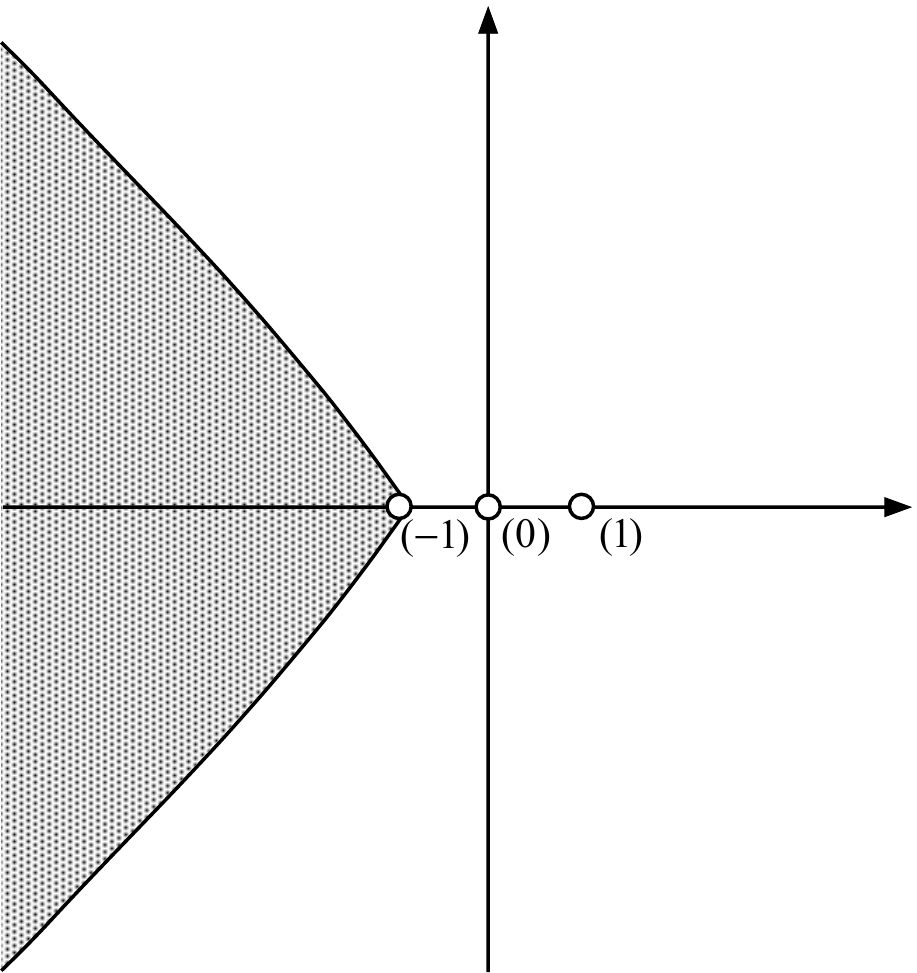}
  \caption{Domain $Z$}
  \label{fig:fig2apcf}
\end{figure}

Let us now record the Airy function approximations using \cite[Thm. 3.4]{Dunster:2020:SEB}. These are valid in an unbounded domain containing the turning point $z=1$. This domain, $Z$ say, contains the points at infinity $z=+\infty$ and $z=\pm i \infty$, but not $z=-\infty$ nor the second turning point $z=-1$. The domain is the one depicted in \cref{fig:fig2apcf}, with boundaries (not included in the domain) the level curves emanating from $z=-1$.

We first define two sequences $\left\{a_{s}\right\} _{s=1}^{\infty}$ and $\left\{\tilde{a}_{s}\right\} _{s=1}^{\infty}$ by $a_{1}=a_{2}=\frac{5}{72}$, $\tilde{a}_{1}=\tilde{a}_{2}=-\frac{7}{72}$, with subsequent terms $a_{s}$ and $\tilde{{a}}_{s}$ ($s=2,3,\cdots $) satisfying the same recursion formula
\begin{equation} \label{3.15}
b_{s+1}=\frac{1}{2}\left(s+1\right) b_{s}+\frac{1}{2}
\sum\limits_{j=1}^{s-1}{b_{j}b_{s-j}}.
\end{equation}
Next let
\begin{equation} \label{3.16}
\mathcal{E}_{s}(z) =\mathrm{E}_{s}(\beta) +
(-1)^{s}a_{s}s^{-1}\xi^{-s},
\end{equation}
and
\begin{equation} \label{3.17}
\tilde{\mathcal{E}}_{s}(z) =\mathrm{E}_{s}(\beta)
+(-1)^{s}\tilde{a}_{s}s^{-1}\xi^{-s}.
\end{equation}

We use the standard notation for Airy functions of complex argument $\mathrm{Ai}_{l}(z):=\mathrm{Ai}(z e^{-2\pi il/3})$ ($l=0,\pm 1$). Then for each nonnegative integer $m$ there exist three solutions of (\ref{3.3}) of the form 
\begin{equation} \label{3.18}
w_{m,l}(u,z) =\mathrm{Ai}_{l}\left( u^{2/3}\zeta \right) 
\mathcal{A}_{2m+2}(u,z) +\mathrm{Ai}_{l}^{\prime }\left(u^{2/3}\zeta\right) \mathcal{B}_{2m+2}(u,z) \ (l=0,\pm 1),
\end{equation}
where
\begin{multline} \label{3.19}
\mathcal{A}_{2m+2}(u,z) =\left( \frac{\zeta }{z^2-1 }\right) ^{1/4}\exp \left\{ \sum\limits_{s=1}^{m}\frac{
\mathcal{\tilde{E}}_{2s}(z) }{u^{2s}}\right\} \cosh \left\{ \sum\limits_{s=0}^{m}\frac{\mathcal{\tilde{E}}_{2s+1}(z) }{u^{2s+1}}\right\}  \\ 
+\frac{1}{2}\left( \frac{\zeta }{z^2-1 }\right) ^{1/4}\tilde{\varepsilon}_{2m+2}(u,z),
\end{multline}
and
\begin{multline} \label{3.20}
\mathcal{B}_{2m+2}(u,z) =\frac{1}{u^{1/3}\left\{\zeta \left(z^2-1\right) \right\}^{1/4}}\exp \left\{ \sum\limits_{s=1}^{m}\frac{
\mathcal{E}_{2s}(z) }{u^{2s}}\right\} \sinh \left\{ \sum\limits_{s=0}^{m}\frac{\mathcal{E}_{2s+1}(z) }{u^{2s+1}}\right\}  \\ 
+\frac{\varepsilon _{2m+2}(u,z) }{2u^{1/3}\left\{ \zeta \left( z^2-1\right) \right\} ^{1/4}}.
\end{multline}

Here it is understood that first sums in (\ref{3.19}) and (\ref{3.20}) are zero if $m=0$. Principal branches are taken for the roots, and both coefficient functions are real for $-1<z< \infty$. The error terms $\tilde{\varepsilon} _{2m+2}(u,z)$ and $\varepsilon _{2m+2}(u,z)$ are $\mathcal{O}(u^{-2m-2})$ uniformly for $z$ lying in the any unbounded closed subset of $Z\setminus\{1\}$. For $z \in Z\setminus\{1\}$ they are bounded by
\begin{multline} \label{3.21}
\left\vert \varepsilon _{2m+2}(u,z) \right\vert \leq \frac{1}{u ^{2m+2}}\exp \left\{ \sum\limits_{s=1}^{2m+1}
{\frac{\Re \left\{\mathcal{E}_{s}(z) \right\}}{u^{s}}}\right\} e_{2m+2,j}(u,z) \left\{ 1+\frac{e_{2m+2,j}(u,z) }{2{u ^{2m+2}}}\right\} ^{2} \\ +\frac{1}{{u ^{2m+2}}}\exp \left\{ \sum\limits_{s=1}^{2m+1}{(-1) ^{s}\frac{\Re \left\{ \mathcal{E}_{s}(z) \right\} }{u^{s}}}\right\} e_{2m+2,k}(u,z) \left\{1+\frac{e_{2m+2,k}(u,z)}{2{u ^{2m+2}}}
\right\} ^{2},
\end{multline}
where for $j,k=0,\pm 1$, $j<k$
\begin{multline} \label{3.22}
e_{n,j}(u,z) ={u ^{n}}\left\vert {\delta 
}_{n,j}(u) \right\vert  
+\omega _{n,j}(u,z) \exp \left\{ {u
^{-1}\varpi _{n,j}(u,z) +u ^{-n}\omega
_{n,j}(u,z) }\right\}  \\
+\gamma _{n}(u,\xi) \exp \left\{ {u
^{-1}\beta _{n}(u,\xi) +u ^{-n}\gamma
_{n}(u,\xi) }\right\}.
\end{multline}
Here $\delta _{n,0}(u)=0$ and
\begin{equation} \label{30.1}
\delta _{n,\pm 1}(u) =\lambda _{\pm 1}
\exp \left\{-2\sum\limits_{s=1}^{n-1}{
\frac{\mathrm{E}_{2s+1}(1)}
{u^{2s+1}}}\right\}-1
=\mathcal{O}(u^{-n}),
\end{equation}
where $\lambda _{\pm 1}(u)$ are the connection coefficients in the relation
\begin{equation} \label{30.2}
\lambda _{-1}(u)W_{-1}(u,\zeta) =iW_{0}(u,\zeta)
+\lambda _{1}(u)W_{1}(u,\zeta),
\end{equation}
where $W_{j}(u,\zeta)$ are the unique solutions of (\ref{3.3}) having the behaviour
\begin{equation}
\lim_{z\rightarrow +\infty}z^{1/2}
e^{u\xi }W_{0}(u,\zeta) =1,
\label{eq14a}
\end{equation}
and
\begin{equation}
\lim_{z\rightarrow \mp i \infty}z^{1/2}
e^{-u\xi }W_{\pm 1}(u,\zeta) =1.
\label{eq14b}
\end{equation}
We shall shortly show that
\begin{equation} \label{30.3}
\lambda _{1}(u)=\lambda _{-1}(u)
=\left( \frac{2e}{u} \right)^{u/2}
\frac{\Gamma\left(\tfrac{1}{2}u+\tfrac{1}{2}\right)}
{\sqrt{2\pi}}.
\end{equation}

Also in (\ref{3.22}) $\omega_{n,j}(u,\beta)$ and $\varpi _{n,j}(u,\beta)$ are given by (\ref{2.17}) and (\ref{2.18}) respectively, with $\bar{\beta}$ replaced by $\beta$, and $(-1)^j$ replaced by 1. The integrals now must be taken along paths that, under the map (\ref{3.7}), are paths in the $z$ plane that avoid $\pm 1$ (i.e. where $\beta$ is unbounded), and on which $\Re(\xi)$ is monotonic. They must also consist of a finite chain of $R_{2}$ arcs. All points in $Z$ must be accessible to at least two of $z=+\infty$ and $z=\pm i \infty$ by such a path. In fact all points in $Z$ in the upper half plane are accessible to $z=+\infty$ and $z= i \infty$, and all points in $Z$ in the lower half plane are accessible to $z=+\infty$ and $z= -i \infty$. Thus in (\ref{3.21}) and (\ref{3.22}) we can take $j=0$, $k=1$ in the former case, and $j=-1$, $k=0$ in the latter.

The bounds break down near $z=1$ since the coefficients $\mathcal{E}_{s}(z)$ and $\tilde{\mathcal{E}}_{s}(z)$ are unbounded there, but crucially these bounds are valid in a deleted neighbourhood of this turning point. Moreover, the coefficient functions $\mathcal{A}_{2m+2}(u,z)$ and $\mathcal{B}_{2m+2}(u,z)$ are analytic there. Hence near, and at $z=1$, one can instead use Cauchy's integral formula
\begin{equation} \label{3.30}
\mathcal{A}_{2m+2}(u,z)=\frac{1}{2\pi i}\oint_{\left| {t-1 } 
\right|=r_{0} } \frac{\mathcal{A}_{2m+2}(u,t)
dt}{t-z} \quad (r_{0}<2),
\end{equation}
and similarly for $\mathcal{B}_{2m+2}(u,z)$. In these integrals the expansions above are used to compute the functions along the boundary (which is away from the turning point), and likewise the error terms can be bounded in a similar way. See \cite[Thm. 4.2]{Dunster:2020:SEB} for details, including error bounds. This method has been shown to be highly accurate and numerically stable for similar applications to Bessel functions \cite{Dunster:2017:COA} and Laguerre polynomials \cite{Dunster:2018:USE}.

Before identifying the asymptotic solutions (\ref{3.18}) with the parabolic cylinder functions, let us prove (\ref{30.3}). Firstly, as $z \rightarrow +\infty$ ($\beta \rightarrow 1$) we have from (\ref{3.4})
\begin{equation} \label{3.23}
\tfrac{2}{3}\zeta^{3/2}=\xi=
\tfrac{1}{2}z^{2}-\tfrac{1}{2}\ln(2z)
-\tfrac{1}{4}
+\mathcal{O}(z^{-2}).
\end{equation}
Hence from (\ref{1.7})
\begin{equation} \label{30.4}
U\left(-\tfrac{1}{2}u,\sqrt {2u} z\right)
\sim \left( \frac{u}{2e} \right)^{u/4}
\frac{e^{-u\xi}}{\left(2uz^2\right)^{1/4}}.
\end{equation}
It follows that by uniqueness of solutions recessive at $z= + \infty$
\begin{equation} \label{30.5}
U\left(-\tfrac{1}{2}u,\sqrt {2u} z\right)
= \left( \frac{u}{2e} \right)^{u/4}
\frac{1}{(2u)^{1/4}}W_{0}(u,z).
\end{equation}

Similarly, for solutions recessive at $z=\pm i \infty$ we find that
\begin{equation} \label{30.6}
U\left(\tfrac{1}{2}u,\pm i\sqrt {2u} z\right)
=\left( \frac{2e}{u} \right)^{u/4}
\frac{e^{\mp (u+1)\pi i/4}}{(2u)^{1/4}}
W_{\pm 1}(u,z).
\end{equation}
So (\ref{30.3}) now follows from inserting (\ref{30.5}) and (\ref{30.6}) into the connection formula \cite[Eq. 12.2.18]{NIST:DLMF}
\begin{multline}
\label{1.8}
\sqrt{2\pi}U(-a,z)=\Gamma\left(\tfrac{1}{2}+a\right)\left\{
e^{\pm \left(\frac{1}{2}a-\frac{1}{4} \right)\pi i}U(a,\pm iz) \right.
\\
+\left. e^{\mp \left(\frac{1}{2}a-\frac{1}{4} \right)\pi i}U\left(a,\mp iz\right)\right\},
\end{multline}
and then comparing this with (\ref{30.2}).

We note in passing that an asymptotic expansion for $U(-\tfrac{1}{2}u,\sqrt {2u} z)$, uniformly valid in an domain that contains the oscillatory interval $[-1+\delta,1-\delta]$ ($\delta \in (0,1)$), can be obtained from (\ref{30.6}) and (\ref{1.8}) by constructing LG expansions for $W_{\pm 1}(u,z)$ similarly to (\ref{2.14}) and  (\ref{2.15}). Such an expansion would be similar to \cite[Eq. (5.11)]{Olver:1959:USE}, but with the advantage that the new one would only involve a single sine (or cosine) term, with the asymptotic series in its argument. Also error bounds would be provided. From this an asymptotic expansion could then be derived for the zeros of $U(-\tfrac{1}{2}u,\sqrt {2u} z)$ (and hence as a special case the Hermite polynomials: see \cite[Eq. 12.7.2]{NIST:DLMF}).  We shall not pursue this here.

Let us now match the asymptotic solutions given by (\ref{3.18}) with the parabolic cylinder functions. From (\ref{3.8}) - (\ref{3.10}), (\ref{3.16}), (\ref{3.17}), (\ref{3.19}), (\ref{3.20}) and (\ref{3.23}) we have as $z \rightarrow \infty$ in the right half plane
\begin{equation} \label{3.24}
\mathcal{A}_{2m+2}(u,z) \sim 
\frac{\zeta^{1/4}}{z^{1/2}}
\cosh \left\{ \sum\limits_{s=0}^{m}\frac{\mathrm{E}_{2s+1}(1)}{u^{2s+1}}+\epsilon_{2m+2}^{(\mathcal{A})}(u)\right\},
\end{equation}
and
\begin{equation} \label{3.25}
\mathcal{B}_{2m+2}(u,z) \sim
\frac{1}{u^{1/3}\zeta^{1/4}z^{1/2}}
\sinh \left\{ \sum\limits_{s=0}^{m}\frac{\mathrm{E}_{2s+1}(1)}{u^{2s+1}}+\epsilon_{2m+2}^{(\mathcal{B})}(u)\right\},
\end{equation}
where $\epsilon_{2m+2}^{(\mathcal{A})}(u)$ and $\epsilon_{2m+2}^{(\mathcal{B})}(u)$ are real constants that are $\mathcal{O}(u^{-2m-2})$ as $u \rightarrow \infty$, and can be bounded with the aid of (\ref{3.21}) (and likewise $\epsilon_{2m+2}^{(U)}(u)$, $\epsilon_{2m+2}^{(\pm 1)}(u)$ and $\epsilon_{2m+2}^{(V)}(u)$ in (\ref{3.28}), (\ref{3.29}) and (\ref{3.34}) below). We omit details since such bounds are cumbersome.

In addition, as $u^{2/3}\zeta \rightarrow \infty$ \cite[Eqs. 9.7.5 and 9.7.6]{NIST:DLMF}
\begin{equation} \label{3.26}
\mathrm{Ai}\left( u^{2/3}\zeta \right) \sim \frac{e^{-u\xi }}{2\pi
^{1/2}u^{1/6}\zeta ^{1/4}}\quad \left(\left\vert \arg\left(u^{2/3}\zeta \right) \right\vert \leq \pi-\delta\right),
\end{equation}
and
\begin{equation} \label{3.27}
\mathrm{Ai}'\left( u^{2/3}\zeta \right) \sim 
-\frac{u^{1/6}\zeta^{1/4}e^{-u\xi }}{2\pi
^{1/2}}\quad \left(\left\vert \arg\left(u^{2/3}\zeta \right) \right\vert \leq \pi-\delta\right).
\end{equation}
Hence from (\ref{1.7}), (\ref{3.18}), (\ref{3.23}), (\ref{3.24}) - (\ref{3.27}) we obtain on matching solutions recessive at $z=+\infty$
\begin{multline} \label{3.28}
U\left(-\tfrac{1}{2}u,\sqrt {2u} z\right)
=\sqrt {\pi }2^{\frac{3}{4}-\frac{1}{4}u}
u^{\frac{1}{4}u-\frac{1}{12}} \\
\times \exp\left\{-\frac{1}{4}u
+\sum\limits_{s=0}^{m}\frac{\mathrm{E}_{2s+1}(1)}{u^{2s+1}}
+\epsilon_{2m+2}^{(U)}(u)\right\}w_{m,0}(u,z),
\end{multline}
and similarly for solutions recessive at $z=\pm i \infty$ we find 
\begin{multline} \label{3.29}
U\left(\tfrac{1}{2}u,\mp i\sqrt {2u} z\right)
=\sqrt {\pi }2^{\frac{3}{4}+\frac{1}{4}u}
u^{-\frac{1}{4}u-\frac{1}{12}} \\
\times \exp\left\{\frac{1}{4}(1 \pm \pi i)u \pm\frac{1}{12}\pi i
-\sum\limits_{s=0}^{m}\frac{\mathrm{E}_{2s+1}(1)}{u^{2s+1}}
+\epsilon_{2m+2}^{(\pm 1)}(u)\right\} w_{m,\pm 1}(u,z).
\end{multline}
Here $\epsilon_{2m+2}^{(\pm 1)}(u)$ are complex conjugates.

Consider next the following numerically satisfactory companion to $U(-a,z)$ for real nonnegative $z$ \cite[Eq. 12.2.20]{NIST:DLMF}
\begin{equation} \label{3.32}
V(-a,z)=(2\pi)^{-1/2}
\left\{
e^{\left(\frac{1}{2}a+\frac{1}{4}\right)\pi i}U(a,iz)
+e^{-\left(\frac{1}{2}a+\frac{1}{4}\right)\pi i}U(a,-iz)
\right\}.
\end{equation}
On using the Airy function connection formula \cite[Eq. 9.2.10]{NIST:DLMF}
\begin{equation} \label{3.33}
\mathrm{Bi}(z)= 
e^{-\frac{1}{6}\pi i}\mathrm{Ai}_{1}(z)
+e^{\frac{1}{6}\pi i}\mathrm{Ai}_{-1}(z),
\end{equation}
we find from (\ref{3.18}), (\ref{3.29}), (\ref{3.32}), (\ref{3.33}) that for $z \in Z$ (and in particular for $0 \leq z < \infty$) 
\begin{multline} \label{3.34}
V\left(-\tfrac{1}{2}u,\sqrt {2u} z\right)
=2^{\frac{1}{4}+\frac{1}{4}u}
u^{-\frac{1}{4}u-\frac{1}{12}}
\exp\left\{\frac{1}{4}u
-\sum\limits_{s=0}^{m}\frac{\mathrm{E}_{2s+1}(1)}{u^{2s+1}}
+\epsilon_{2m+2}^{(V)}(u)
\right\}
\\
\times \left\{\mathrm{Bi}\left( u^{2/3}\zeta \right) 
\mathcal{A}_{2m+2}(u,z) +\mathrm{Bi}^{\prime }\left(u^{2/3}\zeta\right) \mathcal{B}_{2m+2}(u,z)
\right\}
\end{multline}

For extensions covering the left-half plane, and in particular at the turning point $z=-1$ we can use the above approximations along with the connection formulas \cite[Eqs. 12.2.16 and 12.2.19]{NIST:DLMF}
\begin{equation} \label{3.34a}
U(-a,-z)=\mp ie^{\pm i\pi a}U(-a,z)
+\frac{\sqrt{2\pi}}
{\Gamma\left(\tfrac{1}{2}-a\right)}
e^{\pm i\pi(\frac{1}{2}a+\frac{1}{4})}
U\left(a,\pm iz\right),
\end{equation}
and
\begin{equation} \label{3.34b}
V(-a,-z)=\frac{\cos(\pi a)}
{\Gamma\left(\frac{1}{2}+a\right)}U(-a,z)
-\sin(\pi a)V(-a,z).
\end{equation}

For the corresponding uniform asymptotic approximations of the derivatives of the parabolic cylinder functions one can use the above results and appropriate connection formulas given in \cite[Sect. 12.8(i)]{NIST:DLMF}.

\subsection{Inhomogeneous equation}
Here the inhomogeneous version of (\ref{3.1}) that we study is given by
\begin{equation} \label{3.35}
\frac{d^{2}y}{dz^{2}}-\left(\frac{1}{4}z^{2}-a \right)y
=z^{R}.
\end{equation}

Again the fundamental particular solutions are given by (\ref{2.33}), but with $a$ replaced by $-a$. Our main interest is $U_{R}^{(0,2)}(-a,z)$ given by (\ref{2.36}) with $a$ replaced by $-a$, since this is the unique solution that is not exponentially large in a domain containing the whole real $z$ axis. Note from (\ref{2.36}) we have $U_{R}^{(0,2)}(-a,-z)=(-1)^{R}U_{R}^{(0,2)}(-a,z)$. 

Now since $U(-a,z)$ and $U(-a,-z)$ are linearly dependent when $a-\frac{1}{2}=0,1,2,\cdots$, and in particular \cite[Eq. 12.2.15	]{NIST:DLMF}
\begin{equation} \label{3.36}
 U\left(-N-\tfrac{1}{2},-z\right) =
(-1)^{N} U\left(-N-\tfrac{1}{2},z\right) \quad (N=0,1,2,\cdots),
\end{equation}
we see that $U_{R}^{(0,2)}(-a,z)$ is undefined if $a=N+\tfrac{1}{2}$ with $N-R$ being even: if $N-R$ is odd $U_{R}^{(0,2)}(-N-\frac{1}{2},z)$ does exist, since in this case (\ref{2.36}) holds in the limiting case $a \rightarrow -N-\frac{1}{2}$.

Similarly to the previous section we define
\begin{equation} \label{3.37}
w(u,z)=
(2u)^{-\frac{1}{2}R-1}y\left(-\tfrac{1}{2}u,\sqrt{2u} z\right),
\end{equation}
where $y(a,z)$ is any solution of (\ref{3.35}), and then $w(u,z)$ satisfies
\begin{equation} \label{3.38}
\frac{d^{2}w}{dz^{2}}-u^{2}(z^{2}-1)w=z^{R}.
\end{equation}

We then apply \cite[Thm. 4]{Dunster:2020:ASI} to obtain the asymptotic solutions
\begin{equation} \label{3.39}
w_{R}^{(j,k)}(u,z)=\frac{1}{u^{2}}\sum\limits_{s=0}^{n-1} \frac{G_{s,R}(z)}{u^{2s}} +\varepsilon_{n,R}^{(j,k)}(u,z),
\end{equation}
where
\begin{equation} \label{3.40}
G_{0,R}(z)=-z^{R}/\left(z^{2}-1\right),
\end{equation}
and
\begin{equation} \label{3.41}
G_{s+1,R}(z)=G''_{s,R}(z)/\left(z^{2}-1\right)
\quad (s=0,1,2,\cdots ).
\end{equation}
For $n>\frac{1}{4}R-\frac{3}{8}$ the error terms satisfy the bounds (\ref{2.43}) with bars removed, and on the integration path $\Re(\xi)$ varies continuously and is monotonic.

Next from \cite[Thm. 8]{Dunster:2020:ASI} there is a constant $\gamma_{m,R}(u)$ ($m=1,2,3,\cdots$) such that
\begin{equation} \label{3.47}
w_{R}^{(-1,1)}(u,z)=w_{R}^{(0,1)}(u,z)+2\pi e^{-\pi i/6}\gamma_{m,R}(u)w_{m,1}(u,z),
\end{equation}
\begin{equation} \label{3.48}
w_{R}^{(-1,0)}(u,z)=w_{R}^{(-1,1)}(u,z)-2\pi e^{\pi i/6}\gamma_{m,R}(u)w_{m,-1}(u,z),
\end{equation}
and
\begin{equation} \label{3.49}
w_{R}^{(-1,0)}(u,z)=w_{R}^{(0,1)}(u,z)-2\pi i\gamma_{m,R}(u)w_{m,0}(u,z).
\end{equation}
The precise value is given by the following.
\begin{lemma}
\begin{multline} \label{3.50}
\gamma_{m,R}(u)=
2^{-\frac{1}{2}u+R-\frac{1}{2}}
u^{\frac{1}{4}u-\frac{1}{2}R-\frac{13}{12}}
\sqrt {\pi} \exp\left\{-\frac{1}{4}u
+\sum\limits_{s=0}^{m}\frac{\mathrm{E}_{2s+1}(1)}{u^{2s+1}}
+\epsilon_{2m+2}^{(U)}(u)\right\}
\\
\times
\mathbf{F}\left(\tfrac{1}{2}-\tfrac{1}{2}R,
-\tfrac{1}{2}R; \tfrac{1}{4}u-\tfrac{1}{2}R+\tfrac{3}{4};
\tfrac{1}{2}\right).
\end{multline}
\end{lemma}
\begin{remark}
With the aid of Stirling's formula \cite[Eq. 5.11.3]{NIST:DLMF} and (\ref{2.46}) we find that $\gamma_{m,R}(u) \sim 2^{-1/2}u^{-4/3}$ as $u \rightarrow \infty$.
\end{remark}

\begin{proof}

Assume for a moment that $z$ is real, taking the conjugate of (\ref{2.47}), and bearing in mind that in this case $U_{R}^{(0,2)}(a,z)$ and $U(a,z)$ are real, and that the complex conjugate $\overline{U_{R}^{(0,1)}(a,z)}=U_{R}^{(0,3)}(a,z)$, we find that
\begin{equation} \label{3.52}
U_{R}^{(0,2)}(a,z)
=U_{R}^{(0,3)}(a,z)
+\overline{\Lambda_{R}(a)}U(a,z),
\end{equation}
where $\overline{\Lambda_{R}(a)}$ is the complex conjugate of $\Lambda_{R}(a)$. Therefore, on eliminating $U_{R}^{(0,2)}(a,z)$ from (\ref{2.47}) and (\ref{3.52}), and then replacing $a$ by $-a$, we get
\begin{equation} \label{3.53}
U_{R}^{(0,3)}(-a,z)
=U_{R}^{(0,1)}(-a,z)
+2i\Im\left\{\Lambda_{R}(-a)\right\} U(-a,z).
\end{equation}
From (\ref{2.56})
\begin{equation} \label{3.54}
\Im\left\{\Lambda_{R}(-a)\right\}
=-2^{-\frac{1}{2}a+\frac{3}{2}R
-\frac{1}{4}} \pi 
\mathbf{F}\left(\tfrac{1}{2}-\tfrac{1}{2}R,
-\tfrac{1}{2}R;
\tfrac{1}{2}a-\tfrac{1}{2}R+\tfrac{3}{4};\tfrac{1}{2}
\right).
\end{equation}
On recalling that solutions of (\ref{3.35}) and (\ref{3.38}) are related by (\ref{3.37}) we deduce that by matching the unique solution that is not exponentially large at $z=+\infty$ and $z=i\infty$
\begin{equation} \label{3.56}
U_{R}^{(0,1)}\left(-\tfrac{1}{2}u,\sqrt{2u} z\right)
=(2u)^{\frac{1}{2}R+1}w_{R}^{(0,1)}(u,z).
\end{equation}

Similarly
\begin{equation} \label{3.57}
U_{R}^{(0,3)}\left(-\tfrac{1}{2}u,\sqrt{2u} z\right)
=(2u)^{\frac{1}{2}R+1}w_{R}^{(0,-1)}(u,z).
\end{equation}
So from (\ref{3.28}), (\ref{3.49}), (\ref{3.53}) - (\ref{3.57}) we get (\ref{3.50}).
\end{proof}

We can match the asymptotic solutions given by (\ref{3.39}) - (\ref{3.41}) with the corresponding particular solutions which are not exponentially large in the same domains (see (\ref{3.56}) and (\ref{3.57}) below). However our main priority is to obtain asymptotic expansions which are valid at the turning point $z=1$. These involve Scorer functions, defined as follows. Firstly, $\mathrm{Hi}(z)$ is defined by
\begin{equation}
\label{eqHi}
\mathrm{Hi}(z)=\frac{1}{\pi }\int_{0}^{\infty} \exp  \left(-\tfrac{1}{3} t^{3}+z t \right)dt,
\end{equation}
and is the uniquely defined particular solution of the inhomogeneous Airy equation
\begin{equation}
\label{eqHiEq}
\frac{d^{2}w}{dz^{2}}-zw=\frac{1}{\pi},
\end{equation}
having the behaviour
\begin{equation}
\label{HiZ}
\mathrm{Hi}(z)\sim -\frac{1}{\pi z} \quad \left(z\to \infty,\; \left| \arg (-z) \right|\le \tfrac{2}{3}\pi -\delta \right),
\end{equation}
for arbitrary small positive $\delta$. See \cite[Sect. 9.12]{NIST:DLMF} for further properties. Following \cite{Dunster:2020:ASI} we then define
\begin{equation} \label{3.42}
 \mathrm{Wi}^{(-1,1)}(z)=\pi \mathrm{Hi}(z),
\end{equation}
 \begin{equation}  \label{3.43}
 \mathrm{Wi}^{(0,1)}(z)=\pi e^{-2\pi i/3}\mathrm{Hi}\left(ze^{-2\pi i/3}\right),
\end{equation}
and
\begin{equation} \label{3.44}
 \mathrm{Wi}^{(-1,0)}(z)=\pi e^{2\pi i/3}\mathrm{Hi}\left(ze^{2\pi i/3} 
\right).
\end{equation}

These three particular solutions of the inhomogeneous Airy equation equation $w''-zw=1$ are the unique ones which are bounded as $z \rightarrow \infty$ in the sectors $|\arg(-z e^{-2\pi(j+k) i/3})|\le \tfrac{2}{3}\pi -\delta$. For asymptotic expansions with error bounds see \cite[Thm. 5]{Dunster:2020:ASI}. 

Now we are in a position to present expansions for $w^{(j,k)}_{R}(u,z)$ (and hence \\ $U_{R}^{(j,k)}(-a,z)$ from (\ref{3.56}) and (\ref{3.57})) which are valid at the turning point. From \cite[Thm. 7, Eqs (116), (117)]{Dunster:2020:ASI} these read as follows. For any nonnegative integer $m$
\begin{multline} \label{3.58}
 w^{(j,k)}_{R}(u,z)=\gamma_{m,R}(u)\left\{\mathrm{Wi}^{(j,k)}\left(u^{2/3}\zeta  \right)\mathcal{A}_{2m+2}(u,z) \right. \\ + \left. {\mathrm{Wi}^{(j,k)}}^{\prime}\left(u^{2/3}\zeta \right)\mathcal{B}_{2m+2}(u,z) \right\} +\mathcal{G}_{m,R}(u,z) \quad (j,k=0,\pm1, \, j<k),
\end{multline}
where
\begin{multline} \label{3.59}
\mathcal{G}_{m,R}(u,z)=\frac{1}{u^{2}}
\sum\limits_{s=0}^m
\frac{G_{s,R}^{\ast}(z)}{u^{2s}}
\\ 
 -\frac{\gamma_{m}(u)}{2\pi iu^{2/3}}\oint_{\Gamma} {\left\{ {\frac{\zeta(t)}{t^2-1}} \right\}^{1/4}\frac{J_{m}(u,t)dt}{\zeta(t)(t-z)}} +\mathcal{O}\left( \frac{1}{u^{2m+4}} \right),
\end{multline}
and
\begin{multline} \label{3.60}
 J_{m}(u,z)=-\exp \left\{\sum\limits_{s=1}^{m} \frac{\tilde{\mathcal{E}}_{2s}(z)}{u^{2s}} \right\}\cosh \left\{\sum\limits_{s=0}^{m} 
\frac{\tilde{\mathcal{E}}_{2s+1}(z)}{u^{2s+1}}  
\right\}\sum\limits_{k=0}^{m} \frac{(3k)!}{k!\left( 3u^{2}\zeta^{3} \right)^{k}} \\ 
 +\frac{1}{u\zeta^{3/2}}\exp \left\{\sum\limits_{s=1}^m \frac{{\mathcal{E}}_{2s}(z)}{u^{2s}}  \right\}\sinh \left\{\sum\limits_{s=0}^{m} 
\frac{{\mathcal{E}}_{2s+1}(z)}{u^{2s+1}}  \right\}\sum\limits_{k=0}^{m}
\frac{(3k+1)!}{k!\left(3u^{2}\zeta^{3} \right)^{k}}.
\end{multline}

In (\ref{3.58}) $\mathcal{A}_{2m+2}(u,z)$ and $\mathcal{B}_{2m+2}(u,z)$ are the turning point coefficient functions given by (\ref{3.19}) and (\ref{3.20}). In (\ref{3.59}) the coefficients $G_{s,R}^{\ast}(z)$ are the analytic parts at $z=1$ (see \cite[Thm. 9]{Dunster:2020:ASI}) of $G_{s,R}(z)$ (given by (\ref{3.40}) and (\ref{3.41})). The first and third sums in (\ref{3.60}) are taken to be zero if $m=0$. Again the expansions (\ref{3.58}) are valid for $z \in Z$ (\cref{fig:fig2apcf}). For a discussion about error bounds and the computation of $\mathcal{G}_{m,R}(u,z)$ near the turning point see \cite[Sect. 5.1]{Dunster:2020:ASI}

Finally we obtain an approximation for $U_{R}^{(0,2)}(-a,z)$ via the connection formula
\begin{equation} \label{3.64}
U_{R}^{(0,2)}(-a,z)
=\tfrac{1}{2}\left\{U_{R}^{(0,1)}(-a,z)
+U_{R}^{(0,3)}(-a,z)\right\}
+\Re\left\{\Lambda_{R}(-a)\right\}U(-a,z),
\end{equation}
which comes from (\ref{2.47}) and (\ref{3.52}). Now from (\ref{2.56})
\begin{multline} \label{3.65}
\Re\left\{\Lambda_{R}(-a)\right\}
=-\left\{\tan(\pi a)
+(-1)^{R}\sec(\pi a)\right\}
\\
\times 2^{-\frac{1}{2}a+\frac{3}{2}R-\frac{1}{4}} \pi 
\mathbf{F}\left(\tfrac{1}{2}-\tfrac{1}{2}R,
-\tfrac{1}{2}R;
\tfrac{1}{2}a-\tfrac{1}{2}R+\tfrac{3}{4};\tfrac{1}{2}
\right).
\end{multline}
And then the desired expansion comes from plugging this, along with (\ref{3.56}), (\ref{3.57}) and (\ref{3.58}) into (\ref{3.64}).

\section{Weber functions with positive parameter} 
\label{sec4}

We now turn our attention to the Weber equation
\begin{equation} \label{4.1}
\frac{d^{2}y}{dz^{2}}-\left(a-\frac{1}{4}z^{2} \right)y=0,
\end{equation}
where $a$ is positive and large, which has solutions
\begin{equation} \label{4.2}
W_{j}(a,z)
=U\left((-1)^{j}ia,(-i)^{j}ze^{-\pi i/4}\right)
\quad (j=0,1,2,3).
\end{equation}
From (\ref{1.7}) we observe that $W_{j}(a,z)$ is recessive at $z=e^{\pi i/4}i^{j}\infty$. Note also that by uniqueness $W_{2}(a,z)=W_{0}(a,-z)$ and $W_{1}(a,z)=W_{3}(a,-z)$, and so we can primarily focus on $W_{0}(a,z)$ and $W_{3}(a,z)$. The important property of these two is that as $z \rightarrow \infty$ for $j=0$ and $j=3$ from (\ref{1.7}) and (\ref{4.2})
\begin{equation} \label{4.3}
W_{j}(a,z)\sim 
e^{-\tfrac{1}{4}\pi a \pm \tfrac{1}{8}\pi i}
z^{\mp ia-\tfrac{1}{2}}
e^{\pm \tfrac{1}{4}iz^{2}}
\quad \left(
\left\vert \arg\left( ze^{\mp\pi i/4}
\right) \right\vert
\leq \tfrac{3}{4}\pi-\delta \right),
\end{equation}
where the upper signs are taken for $j=0$ and lower signs for $j=3$. 

For real argument $x$ the solutions $W(a,\pm x)$ are of most interest, which are defined from equating real and imaginary parts of the equation (see \cite{Olver:1959:USE})
\begin{equation} \label{4.3a}
k^{-1/2}(2a)W(a,x)
+ik^{1/2}(2a)W(a,-x) \\
=2^{1/2}e^{\pi a/4}
e^{i\rho(2a)} W_{0}(a,x),
\end{equation}
where
\begin{equation} \label{4.5}
k(u)=\sqrt {1+e^{\pi u}} -e^{\pi u/2},
\end{equation}
\begin{equation} \label{4.6}
1/k(u)=\sqrt {1+e^{\pi u}} +e^{\pi u/2},
\end{equation}
\begin{equation} \label{4.7}
\rho(u) = \tfrac{1}{2}\phi_{2}(u)+\tfrac{1}{8}\pi,
\end{equation}
and
\begin{equation} \label{4.8}
\phi_{2}(u) =\arg \left\{\Gamma\left(\tfrac{1}
{2}+\tfrac{1}{2}iu
\right) \right\}.
\end{equation}
In (\ref{4.8}) the branch of $\arg$ is taken to be zero when $u=0$ and then defined by continuity for $u>0$. Note for real $u$
\begin{equation} \label{4.9}
e^{2i\phi_{2}(u)} =
\frac{\Gamma\left(\tfrac{1}{2}+\tfrac{1}{2}iu \right)}
{\Gamma\left(\tfrac{1}{2}-\tfrac{1}{2}iu \right)}.
\end{equation}

For complex argument $z$ we also have the relation \cite[Eq. 12.14.4]{NIST:DLMF}
\begin{equation} \label{4.4}
W(a,z)=\sqrt{\frac{k(2a)}{2}} e^{\pi a/4}\left\{ e^{i\rho(2a)
}W_{0}(a,z)+e^{-i\rho(2a) }W_{3}(a,z) \right\}.
\end{equation}
We shall obtain LG and Airy expansions for the complex valued solutions $W_{j}(a,z)$, and use these to obtain the desired expansions for $W(a,x)$ for $-\infty < x < \infty$.

\begin{figure}[htbp]
  \centering
  \includegraphics[width=\textwidth,keepaspectratio]{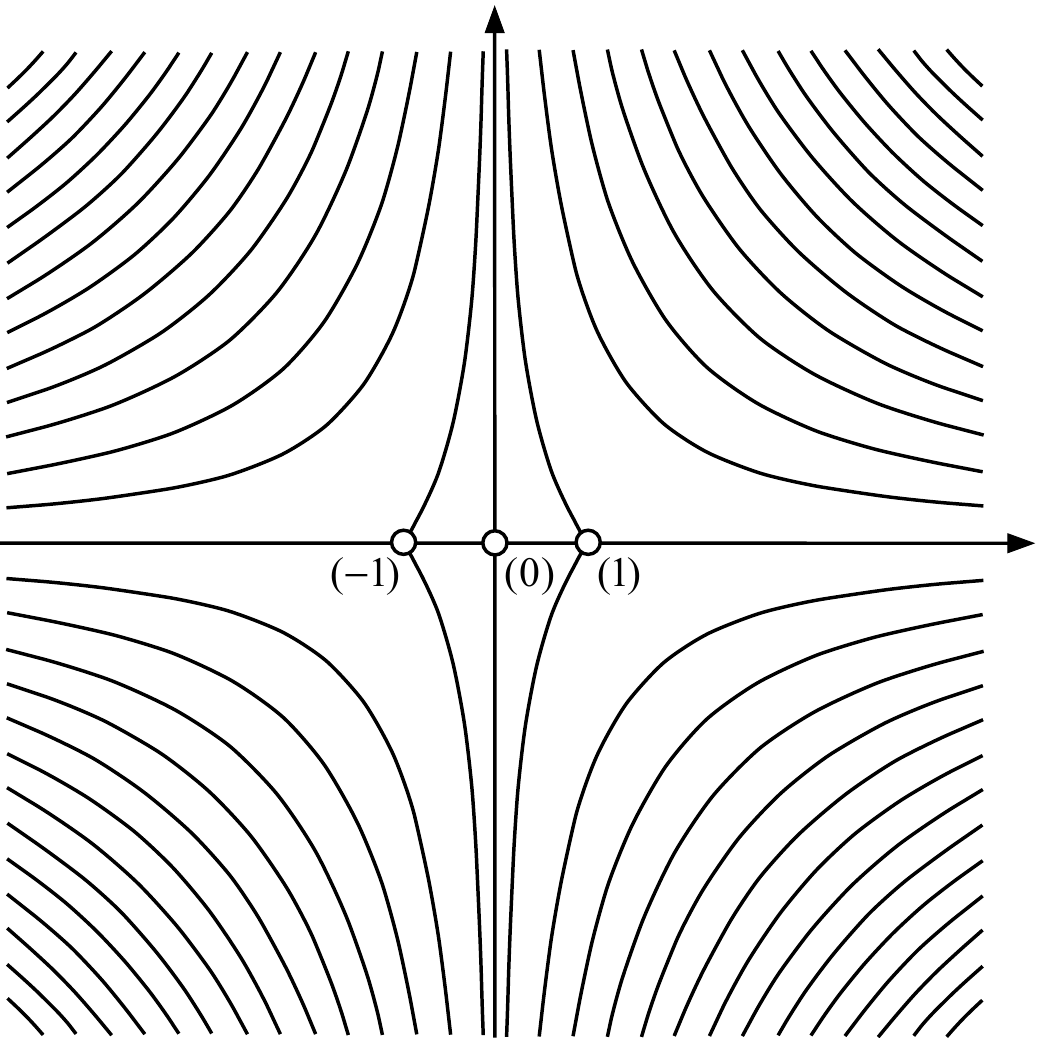}
  \caption{Level curves $\Im(\xi)=\mathrm{constant}$}
  \label{fig:fig3pcf}
\end{figure}

Again with $a=\tfrac{1}{2}u$, $z$ replaced by $\sqrt{2u} z$, and $w(u,z)=u^{-1}y(\tfrac{1}{2}u,\sqrt {2u} z)$. we recast (\ref{4.1}) into
\begin{equation} \label{4.10}
\frac{d^{2}w}{dz^{2}}-u^{2}(1-z^{2})w=0.
\end{equation}
The analysis is similar to \cref{sec3} since the equation has real turning points at $z = \pm 1$. Here the role of $\xi$ is played by $i\xi$, and the level curves $\Re(i\xi)=\Im(\xi)=\mathrm{constant}$ are depicted in \cref{fig:fig3pcf}. These of course are orthogonal to those from \cref{sec3} (\cref{fig:fig2pcf}). The Airy variable $\zeta$ is replaced by $-\zeta$.

We follow the same steps as \cref{sec3}, using $f(z)=1-z^2$ instead of the negative of this function. As a result we find $\mathrm{E}_{s}(\beta)$ are replaced by $(-i)^{s}\mathrm{E}_{s}(\beta)$ (equivalently $i^{s}\bar{\mathrm{E}}_{s}(\beta)$), $\mathcal{E}_{s}(z)$ replaced by $(-i)^{s}\mathcal{E}_{s}(z)$, $\tilde{\mathcal{E}}_{s}(z)$ replaced by $(-i)^{s}\tilde{\mathcal{E}}_{s}(z)$. The solutions valid at the turning point $z=1$ are then given by
\begin{equation} \label{4.11}
w_{m,l}^{(W)}(u,z) =\mathrm{Ai}_{l}\left( -u^{2/3}\zeta\right) 
\mathcal{A}_{2m+2}^{(W)}(u,z) +\mathrm{Ai}_{l}^{\prime }\left(-u^{2/3}\zeta\right) \mathcal{B}_{2m+2}^{(W)}(u,z)
\ (l=0,\pm 1),
\end{equation}
where for nonnegative $m$ $\mathcal{A}_{2m+2}^{(W)}(u,z)$ and $ \mathcal{B}_{2m+2}^{(W)}(u,z)$ are again given by (\ref{3.19}) and (\ref{3.20}), but with $\mathcal{E}_{s}(z)$ and $\tilde{\mathcal{E}}_{s}(z)$ being replaced by $(-i)^{s}\mathcal{E}_{s}(z)$ and $(-i)^{s}\tilde{\mathcal{E}}_{s}(z)$, respectively. Also the error terms, $\varepsilon_{2m+2}^{(W)}(u,z)$ and $\tilde{\varepsilon}_{2m+2}^{(W)}(u,z)$, say, satisfy the bounds (\ref{3.21}) and (\ref{3.22}) with $\mathcal{E}_{s}(z)$ and $\tilde{\mathcal{E}}_{s}(z)$ replaced as above. In these bounds $\omega_{n,j}(u,\beta)$ and $\varpi _{n,j}(u,\beta)$ are the same as described in \cref{sec3}, that is, given by (\ref{2.17}) and (\ref{2.18}) respectively, with $\bar{\beta}$ replaced by $\beta$, and $(-1)^j$ replaced by 1. The only difference is that on the paths of integration $\Im(\xi)$ must be monotonic (and again avoid $z=\pm 1$). 

\begin{figure}[htbp]
  \centering
  \includegraphics[width=\textwidth,keepaspectratio]{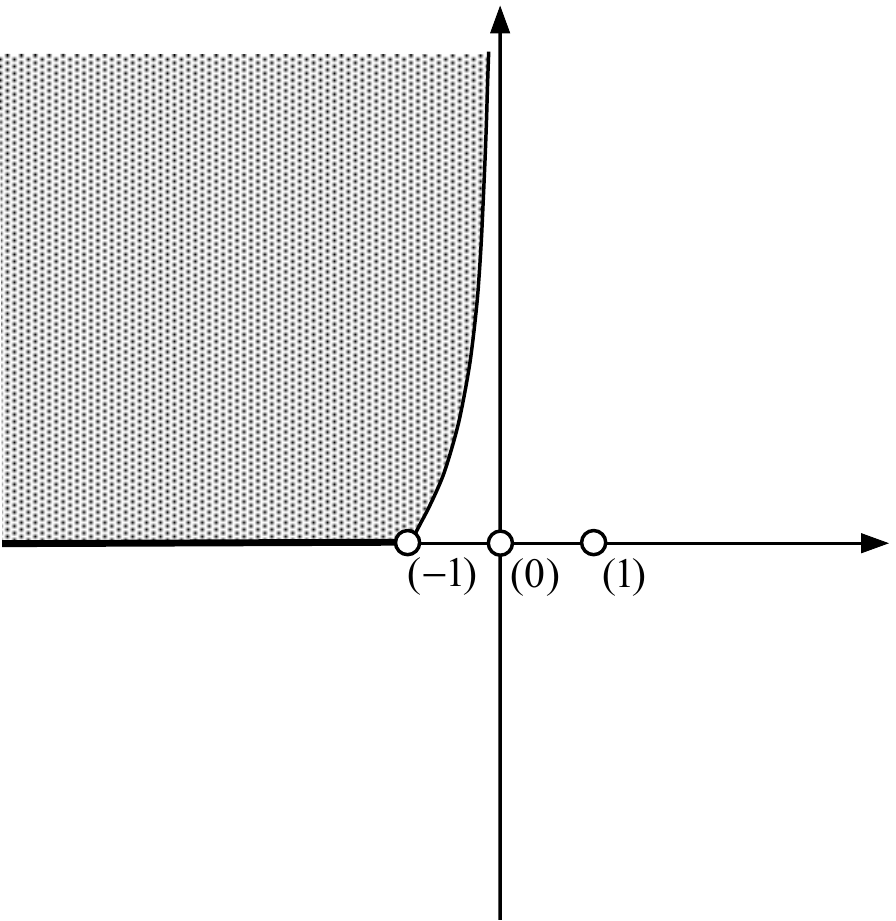}
  \setlength{\belowcaptionskip}{10pt}
  \caption{Domain $\breve{Z}$}
  \label{fig:fig3bpcf}
\end{figure}

We consider $W_{j}(a,z)$ for $j=0,2,3$. For $j=0$ the equivalent path in the $z$ plane has its end point at $z=e^{\pi i/4}\infty$, for $j=2$ the end point is at $z=e^{-3\pi i/4}\infty$, and for $j=3$ the end point is at $z=e^{-\pi i/4}\infty$. All points in the common domain of validity, which we call $\breve{Z}$, must be accessible to at least two of these points at infinity by such a path on which the continuous branch of $\Im(\xi)$ is monotonic. The collection of all such points is denoted by $\breve{Z}$ as shown in \cref{fig:fig3bpcf}. In this figure the boundaries are the level curve $\Im(\xi) = \text{constant}$ and the interval $-\infty < z \leq -1$, the latter being the branch cut for $\zeta$. Note these boundaries do not lie in $\breve{Z}$ (thus it is open). Note that $z=1$ is excluded from the error bounds, since the integrals diverge there. But near this turning point we again can use Cauchy integrals of the form (\ref{3.30}) for both computing the coefficient functions as well as obtaining error bounds.

Since the coefficient functions $\mathcal{A}_{2m+2}^{(W)}(u,z)$ and $\mathcal{B}_{2m+2}^{(W)}(u,z)$ are entire and real on the real axis, we can use the Schwarz reflection principle to immediately extend the domain of asymptotic validity of the asymptotic expansions and their error bounds to the whole $z$ plane, with the exception of points on the cut $(- \infty, -1]$.

We now identify recessive solutions of (\ref{4.10}), and as a result we find for $j=0,3$ that
\begin{multline} \label{4.21}
W_{j}\left(\tfrac{1}{2}u,\sqrt {2u} z\right)
=2^{3/4}\sqrt{\pi}u^{-1/12}
 e^{\pm i (\chi_{m}(u)
+\frac{1}{24}\pi)} \\
\times \exp\left\{-\tfrac{1}{8}\pi u+\epsilon_{2m+2}^{(W \mp)}(u)\right\} 
w_{m,\mp 1}^{(W)}(u,z),
\end{multline}
with upper signs for $j=0$ and lower signs for $j=3$, and $\epsilon_{2m+2}^{(W \pm)}(u) = \mathcal{O}(u^{-2m-2})$. Here we have introduced the parameter
\begin{equation} \label{4.22}
\chi_{m}(u)=
\frac{u}{4}\ln\left(\frac{2e}{u}\right)
+  \sum\limits_{s=0}^{m}
(-1)^{s}\frac{\mathrm{E}_{2s+1}(1)}
{u^{2s+1}}.
\end{equation}

For $j=2$ we likewise find
\begin{multline} \label{4.23}
W_{2}\left(\tfrac{1}{2}u,\sqrt {2u} x\right)
=2^{3/4}\sqrt{\pi}u^{-1/12}
e^{i(\chi_{m}(u)-\frac{1}{8}\pi)}
\\
\times \exp\left\{\tfrac{3}{8}\pi u+\epsilon_{2m+2}^{(W-)}(u)\right\} 
w_{m,0}^{(W)}(u,z).
\end{multline}

Uniform asymptotic approximations for the derivatives of these functions, valid in the same domain, can be derived using (\ref{4.2}), (\ref{4.21}), (\ref{4.23}) and \cite[Sect. 12.8(i)]{NIST:DLMF}.

Let us now record LG expansions, which will primarily be used to relate $\chi_{m}(u)$ to $\phi_{2}(u)$. Similarly to (\ref{2.20}) we find that
\begin{multline} \label{4.24}
W_{j}\left(\tfrac{1}{2}u,\sqrt {2u} z\right)
=\frac{e^{\pm i (\chi_{m}(u)
-\frac{1}{8}\pi)}}
{(2u)^{1/4}e^{\pi u/8}\left(1-z^{2}\right)^{1/4}}
\exp \left\{\sum\limits_{s=1}^{m}(-1)^{s}
\frac{\mathrm{E}_{2s}(\beta)}
{u^{2s}} \right\} 
\\
\times \exp \left\{
\pm i u \xi
\mp i\sum\limits_{s=0}^{m}(-1)^{s}
\frac{\mathrm{E}_{2s+1}(\beta)}
{u^{2s+1}}\right\} 
\left\{1+\eta_{2m+2,j}^{(W)}(u,z) \right\},
\end{multline}
where upper signs for $j=0$ and lower sign for $j=3$. The error terms $\eta_{2m+2,j}^{(W)}(u,z)$ are bounded by (\ref{2.16}) with $\bar{\beta}$ replaced by $\beta$, and $\omega_{n,j}(u,\beta)$ and $\varpi _{n,j}(u,\beta)$ given as described earlier in this section.

\begin{figure}[htbp]
  \centering
  \includegraphics[width=\textwidth,keepaspectratio]{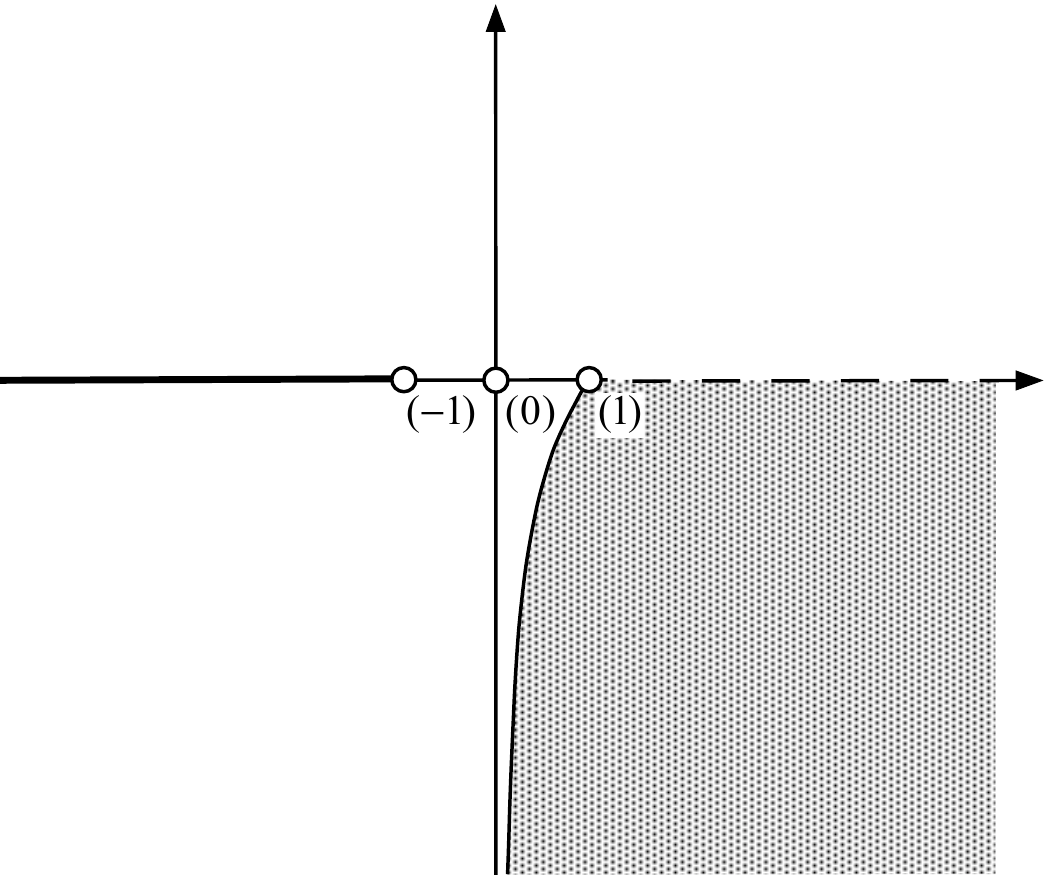}
  \caption{Domain $\breve{Z}_{0}$}
  \label{fig:fig3apcf}
\end{figure}

The domain of validity of (\ref{4.24}) for $W_{0}(\frac{1}{2}u,\sqrt {2u} z)$, $\breve{Z}_{0}$, is shown in \cref{fig:fig3apcf}, with the corresponding one for $W_{0}(\frac{1}{2}u,\sqrt {2u} z)$ being the conjugate of this. Here we take the branch cuts for $\xi$ to be the lines $-\infty < z \leq -1$ (solid line) and $1\leq z < \infty$ (dashed line): the latter can be crossed into the lower half plane without violating the monotonicity condition on $\Im (\xi)$, but we are content to use the unshaded region as shown. The other boundary, the level curve in the fourth quadrant emanating from $z=1$, cannot be crossed: indeed all points on this curve must be excluded from the domain, as well as those on the aforementioned cut associated with $z=-1$.

Now both LG expansions are valid at $z=0$, so from (\ref{4.24}) and \cite[Eq. 12.2.6]{NIST:DLMF}, and recalling that $\bar{\mathrm{E}}_{2s+1}(0)=0$, we have
\begin{multline} \label{4.25}
\frac{W_{0}\left(\tfrac{1}{2}u,0\right)}{W_{3}
\left(\tfrac{1}{2}u,0\right)}
= \frac{U\left(\tfrac{1}{2}iu,0\right)}
{U\left(-\tfrac{1}{2}iu,0\right)}
=\frac{2^{-iu/2}\Gamma\left(
\tfrac{3}{4}-\tfrac{1}{4}iu\right)}
{\Gamma\left(\tfrac{3}{4}+\tfrac{1}{4}iu\right)}
\\
=e^{ i (2\chi_{m}(u) -\frac{1}{4}\pi)}\left\{
\frac{1+\eta_{2m+2,0}^{(W)}(u,0)}
{1+\eta_{2m+2,3}^{(W)}(u,0)}\right\}.
\end{multline}

Next using \cite[Eqs. 5.5.3 and 5.5.5]{NIST:DLMF} one can show with some algebra and simplification that
\begin{equation} \label{4.26}
\frac{2^{-iu/2}\Gamma\left(
\tfrac{3}{4}-\tfrac{1}{4}iu\right)}
{\Gamma\left(\tfrac{3}{4}+\tfrac{1}{4}iu\right)}
=e^{-\pi i/4}
\left\{\frac{\tanh\left(\tfrac{1}{4}\pi (u+i)\right)
\Gamma\left(\tfrac{1}{2}-\tfrac{1}{2}iu\right)}
{\Gamma\left(
\tfrac{1}{2}+\tfrac{1}{2}iu\right)}\right\}^{1/2}.
\end{equation}
Now let
\begin{equation} \label{4.27}
\varepsilon_{m}(u)= \tfrac{1}{2}\phi_{2}(u)+\chi_{m}(u),
\end{equation}
where $\phi_{2}(u)$ is given by (\ref{4.8}). Then from (\ref{4.9}), (\ref{4.25}), (\ref{4.26})
\begin{equation} \label{4.28}
e^{i\varepsilon_{m}(u)}
=\left\{\tanh\left(
\tfrac{1}{4}\pi (u+i)\right)\right\}^{1/4}
\left\{\frac{1+\eta_{2m+2,3}^{(W)}(u,0)}
{1+\eta_{2m+2,0}^{(W)}(u,0)}\right\}^{1/2}
=1+\mathcal{O}(u^{-2m-2}),
\end{equation}
and hence $\varepsilon_{m}(u)=\mathcal{O}(u^{-2m-2})$. Thus from (\ref{4.27}) we arrive at the relation we were aiming for, namely $\chi_{m}(u)= -\tfrac{1}{2}\phi_{2}(u)+\mathcal{O}(u^{-2m-2})$.

We now proceed to obtain our main results, for real $z=x$. From (\ref{4.3a}) we have
\begin{multline} \label{4.29}
k^{-1/2}(u)W\left(\tfrac{1}{2}u,\sqrt {2u} x\right)
+ik^{1/2}(u)W\left(\tfrac{1}{2}u,-\sqrt {2u} x\right) \\
=2^{1/2}e^{\pi u/8}
e^{i\left(\frac{1}{2}\phi_{2}(u)+\frac{1}{8}\pi\right)}
W_{0}\left(\tfrac{1}{2}u,\sqrt {2u}x \right).
\end{multline}
Hence equating real parts, defining $\epsilon_{R,2m+2}(u)=\Re \{\epsilon_{2m+2}^{(W-)}(u)\}$ and $\epsilon_{I,2m+2}(u)=\Im \{\epsilon_{2m+2}^{(W-)}(u)\}$ (both being $\mathcal{O}(u^{-2m-2})$), and using (\ref{4.21}) leads to
\begin{multline} \label{4.30}
W\left(\tfrac{1}{2}u,\sqrt {2u} x\right)
=2^{5/4}\sqrt{\pi k(u)}u^{-1/12} e^{\epsilon_{R,2m+2}(u)}
 \\  \times
\Re\left\{e^{i\left(\frac{1}{2}\phi_{2}(u)+\chi_{m}(u)
+\frac{1}{6}\pi +\epsilon_{I,2m+2}(u)\right)}
w_{m,-1}^{(W)}(u,x) \right\}.
\end{multline}
Thus from (\ref{4.11}), (\ref{4.27}), (\ref{4.28}) and \cite[Eq. 9.2.11]{NIST:DLMF} we arrive at the first of our main asymptotic expansions
\begin{multline} \label{4.31}
W\left(\tfrac{1}{2}u,\sqrt {2u} x\right)
=2^{1/4}\sqrt{\pi k(u)}u^{-1/12} e^{\epsilon_{R,2m+2}(u)}
 \left[
\left\{\mathrm{Bi}\left( -u^{2/3}\zeta\right)
\mathcal{A}_{2m+2}^{(W)}(u,x) \right. \right.\\
\left.  +\mathrm{Bi}^{\prime}\left(-u^{2/3}\zeta\right) \mathcal{B}_{2m+2}^{(W)}(u,x)\right\}
\cos\left(\varepsilon_{m}(u)+\epsilon_{I,2m+2}(u)\right) \\
\left. -\sin\left(\varepsilon_{m}(u)+\epsilon_{I,2m+2}(u)\right)
w_{m,0}^{(W)}(u,x) \right],
\end{multline}
where the error terms associated with the coefficient functions $\mathcal{A}_{2m+2}^{(W)}(u,x)$ and \\ $\mathcal{B}_{2m+2}^{(W)}(u,x)$ are $\mathcal{O}(u^{-2m-2})$ uniformly for $-1+\delta \leq x < \infty$ ($\delta >0$). Also observe that $\cos(\varepsilon_{m}(u)+\epsilon_{I,2m+2}(u))$ = $1+\mathcal{O}(u^{-4m-4})$ and 
\begin{multline} \label{4.31a}
\sin(\varepsilon_{m}(u)+\epsilon_{I,2m+2}(u))w_{m,0}^{(W)}(u,x)
= \mathcal{O}\left(u^{-2m-2}\right)
\mathrm{env} \,
\mathrm{Ai}\left( -u^{2/3}\zeta\right) \\
= \mathcal{O}\left(u^{-2m-2}\right)
\mathrm{env} \, \mathrm{Bi}\left( -u^{2/3}\zeta\right),
\end{multline}
as $u \rightarrow \infty$, where the envelope ($\mathrm{env}$) of the Airy functions are defined by \cite[Eqs. 2.8.20 and 2.8.21]{NIST:DLMF}. 

Thus in summary
\begin{multline} \label{4.31b}
W\left(\tfrac{1}{2}u,\sqrt {2u} x\right)
=2^{1/4}\sqrt{\pi k(u)}u^{-1/12}
 \left[
\left\{\mathrm{Bi}\left( -u^{2/3}\zeta\right)
\mathcal{A}_{2m+2}^{(W)}(u,x) \right. \right.\\
\left. \left. +\mathrm{Bi}^{\prime}\left(-u^{2/3}\zeta\right) \mathcal{B}_{2m+2}^{(W)}(u,x)\right\}
+\mathcal{O}\left(u^{-2m-2}\right)
\mathrm{env} \, \mathrm{Bi}\left( -u^{2/3}\zeta\right)
  \right].
\end{multline}

If we equate imaginary parts of (\ref{4.29}) we are unable to obtain a satisfactory approximation for the numerically satisfactory companion function $W(\tfrac{1}{2}u,-\sqrt {2u}x)$, since an error term dominates the approximating Airy function $\mathrm{Ai}$ and its derivative. We overcome this problem by simply replacing $x$ by $-x$ in (\ref{4.29}) and using $W_{0}(a,-x)=W_{2}(a,x)$. Consequently we have
\begin{multline} \label{4.32}
k^{-1/2}(u)W\left(\tfrac{1}{2}u,-\sqrt {2u}x\right)
+ik^{1/2}(u)W\left(\tfrac{1}{2}u,\sqrt {2u} x\right) \\
=2^{1/2}e^{\pi u/8}
e^{i\left(\frac{1}{2}\phi_{2}(u)+\frac{1}{8}\pi\right)}
W_{2}\left(\tfrac{1}{2}u,\sqrt {2u} x\right).
\end{multline}
Hence from (\ref{4.23}) and by equating real parts of (\ref{4.32}) we obtain the desired expansion
\begin{multline} \label{4.33}
W\left(\tfrac{1}{2}u,-\sqrt {2u}x\right)
=2^{5/4}\sqrt{\pi k(u)}u^{-1/12}e^{\pi u/2}
\left[\mathrm{Ai}\left(-u^{2/3}\zeta\right)
\mathcal{A}_{2m+2}^{(W)}(u,x) \right.
\\   \left.
 +\mathrm{Ai}^{\prime}\left(-u^{2/3}\zeta\right) \mathcal{B}_{2m+2}^{(W)}(u,x)
+\mathcal{O}\left(u^{-2m-2}\right)
\mathrm{env} \, \mathrm{Ai}\left( -u^{2/3}\zeta\right)\right],
\end{multline}
which like (\ref{4.31}) is uniformly valid for $-1+\delta \leq x < \infty$ ($\delta >0$). Of course for both these approximations it suffices to use them for $0 \leq x < \infty$.

\subsection{Inhomogeneous equation}
The inhomogeneous version of (\ref{4.1}) is
\begin{equation} \label{4.34}
\frac{d^{2}y}{dz^{2}}-\left(a-\frac{1}{4}z^{2} \right)y
=z^{R}.
\end{equation}
For $j \in \{0,1,2,3\}$, $k \in \{1,2,3\}$ with $j<k$ solutions of (\ref{4.34}) are furnished by
\begin{multline} \label{4.35}
 W_{R}^{(j,k)}(a,z)
 =\frac{1}
{\mathscr{W}\{W_{j}(a,z),W_{k}(a,z)\}}
\left[ W_{k}(a,z)
 \int_{e^{\pi i/4}i^{j}\infty}^z \right.
 t^{R}W_{j}(a,t)dt \\
\left. -W_{j}(a,z)
 \int_{e^{\pi i/4}i^{k}\infty}^z 
 t^{R}W_{k}(a,t)dt
 \right].
\end{multline}

The solution $W_{R}^{(0,3)}(a,z)$ is the one we are most interested in, because it is the unique solution that is non-oscillatory on the positive real axis, and is $\mathcal{O}(z^{R-2})$ as $z \rightarrow \infty$ in the right half plane, whereas all other solutions are exponentially large in part or all of the right half plane. The function is also real-valued on the real axis. A solution that is numerically satisfactory in the left half plane is $W_{R}^{(1,2)}(a,z)=(-1)^R W_{R}^{(0,3)}(a,-z)$.

Now from (\ref{2.36a}) and (\ref{4.2}) 
\begin{equation} \label{4.35b}
\mathscr{W}\{W_{0}(a,z),W_{3}(a,z)\}
=-ie^{-\pi a/2},
\end{equation}
and hence from (\ref{4.35})
\begin{multline} \label{4.35a}
 W_{R}^{(0,3)}(a,z)
 =ie^{\pi a/2}\left[
W_{3}(a,z)
 \int_{e^{\pi i/4}\infty}^z \right.
 t^{R}W_{0}(a,t)dt \\
\left. -W_{0}(a,z)
 \int_{e^{-\pi i/4}\infty}^z 
 t^{R}W_{3}(a,t)dt
 \right].
\end{multline}

On comparing the differential equations (\ref{3.1}) and (\ref{4.34}), and identifying solutions that are not exponentially large in the same regions of the $z$ plane, we have the following relations for two more solutions that we shall use
\begin{equation} \label{4.37}
W_{R}^{(2,3)}(a,z)
=ie^{-3 R\pi i/4}
U_{R}^{(0,1)}\left(ia,ze^{3\pi i/4}\right),
\end{equation}
and
\begin{equation} \label{4.38}
W_{R}^{(0,2)}(a,z)
=ie^{-3R\pi i/4}
U_{R}^{(0,2)}\left(ia,ze^{3\pi i/4}\right),
\end{equation}
where $U_{R}^{(0,2)}(a,z)$ is given by (\ref{2.36}) and $U_{R}^{(0,1)}(a,z)$ is given by (\ref{2.36b}).

In order to derive a connection relation between these two functions we use (\ref{2.47}) to get
\begin{equation} \label{4.39}
ie^{-3 R\pi i/4}U_{R}^{(0,2)}(ia,z)
=ie^{-3 R\pi i/4}U_{R}^{(0,1)}(ia,z)
+i\Lambda_{R}(ia)e^{-3 R\pi i/4}U(ia,z),
\end{equation}
and hence on replacing $z$ by $ze^{3\pi i/4}$ and using (\ref{4.2}), (\ref{4.37}) and (\ref{4.38}) we arrive at
\begin{equation} \label{4.40}
W_{R}^{(0,2)}(a,z)
=W_{R}^{(2,3)}(a,z)
+ie^{-3 R\pi i/4}\Lambda_{R}(ia)
W_{2}(a,z).
\end{equation}

With $w(u,z)=
(2u)^{-\frac{1}{2}R-1}y(\tfrac{1}{2}u,\sqrt{2u} z)$ we transform (\ref{4.34}) to
\begin{equation} \label{4.42}
\frac{d^{2}w}{dz^{2}}-u^{2}(z^{2}-1)w=z^{R}.
\end{equation}
From \cite{Dunster:2020:ASI} asymptotic solutions valid at $z=1$ are given by
\begin{multline} \label{4.50}
 \breve{w}_{R}^{(j,k)}(u,z)=\gamma_{m,R}^{(W)}(u)\left\{\mathrm{Wi}^{(j,k)}\left(-u^{2/3}\zeta  \right)\mathcal{A}_{2m+2}^{(W)}(u,z) \right. \\ + \left. {\mathrm{Wi}^{(j,k)}}^{\prime}\left(-u^{2/3}\zeta \right)\mathcal{B}_{2m+2}^{(W)}(u,z) \right\} +\mathcal{G}_{m,R}^{(W)}(u,z) \quad (j,k=0,\pm1, \, j<k).
\end{multline}
Here
\begin{multline} \label{4.47}
\gamma_{m,R}^{(W)}(u)=
2^{R-1-\frac{1}{4}iu}
u^{-\frac{1}{2}R-\frac{13}{12}}
e^{\pi u/8}
e^{i(\chi_{m}(u)-\frac{1}{4}\pi R
+\frac{1}{8}\pi )}
\Gamma  \left(\tfrac{1}{2}+\tfrac{1}{2}iu \right)  \\
\times
\mathbf{F}\left(\tfrac{1}{2}-\tfrac{1}{2}R,
-\tfrac{1}{2}R;\tfrac{3}{4}
-\tfrac{1}{2}R+\tfrac{1}{4}iu;\tfrac{1}{2} \right)
\left\{1+\mathcal{O}\left(u^{-2m-2}
\right)\right\},
\end{multline}
$\mathcal{G}_{m,R}^{(W)}(u,z)$ is given by (\ref{3.59}) and (\ref{3.60}) with $\gamma_{m,R}(u)$ replaced by $\gamma_{m,R}^{(W)}(u)$, $G_{s,R}^{\ast}(z)$ are replaced by $(-1)^{s+1}G_{s,R}^{\ast}(z)$, $\zeta$ replaced by $-\zeta$, $t^2-1$ replaced by $1-t^2$, and $\mathcal{E}_{s}(z)$ and $\tilde{\mathcal{E}}_{s}(z)$ are replaced by $(-i)^{s}\mathcal{E}_{s}(z)$ and $(-i)^{s}\tilde{\mathcal{E}}_{s}(z)$, respectively. All error terms are $\mathcal{O}(u^{-2m-2})$ in whole $z$ plane except for points on the cut $(-\infty,-1]$.

The connection coefficient $\gamma_{m,R}^{(W)}(u)$ was derived as follows. Firstly on identifying solutions that are not exponentially large at $z=\infty$ in the first and third quadrants we have
\begin{equation} \label{4.44}
W_{R}^{(0,2)}\left(\tfrac{1}{2}u,\sqrt{2u} z\right)=
(2u)^{\frac{1}{2}R+1}\breve{w}_{R}^{(-1,0)}(u,z),
\end{equation}
and similarly for those that are not exponentially large at $z=\infty$ in the third and fourth quadrants we have
\begin{equation} \label{4.45}
W_{R}^{(2,3)}\left(\tfrac{1}{2}u,\sqrt{2u} z\right)=
(2u)^{\frac{1}{2}R+1}\breve{w}_{R}^{(0,1)}(u,z).
\end{equation}
Now from (\ref{4.40}), (\ref{4.44}) and (\ref{4.45})
\begin{equation} \label{4.46}
\breve{w}_{R}^{(-1,0)}(u,z)
=\breve{w}_{R}^{(0,1)}(u,z)
+ie^{-3 R\pi i/4}
(2u)^{-\frac{1}{2}R-1}
\Lambda_{R}\left(\tfrac{1}{2}iu\right)
W_{2}\left(\tfrac{1}{2}u,z\right).
\end{equation}
So from (\ref{2.55}), (\ref{3.49}) (with $w$ and $\gamma_{m,R}$ replaced by $\breve{w}$ and $\gamma_{m,R}^{(W)}$, respectively), (\ref{4.23}), and (\ref{4.46}) we get (\ref{4.47}).

We note the uniform asymptotic expansion
\begin{equation} \label{4.48}
W_{R}^{(0,3)}\left(\tfrac{1}{2}u,\sqrt{2u} z\right)
=(2u)^{\frac{1}{2}R+1}\breve{w}_{R}^{(-1,1)}(u,z).
\end{equation}
where $\breve{w}_{R}^{(-1,1)}(u,z)$ is given by (\ref{4.50}) with $\mathrm{Wi}^{(-1,1)}(-u^{2/3}\zeta)=\pi \mathrm{Hi}(-u^{2/3}\zeta)$.

\section{Weber functions with negative parameter} 
\label{sec5}
We finally consider Weber's equation in the form
\begin{equation} \label{5.1}
\frac{d^{2}y}{dz^{2}}+\left(\frac{1}{4}z^{2}+a \right)y=0,
\end{equation}
with complex valued solutions $W_{j}(-a,\pm z)$ ($j=0,3$), real solutions $W(-a,\pm x)$ ($-\infty < x < \infty$), and where again $a$ is positive and large.

As before let $a=\tfrac{1}{2}u$, $z$ be replaced by $\sqrt{2u} z$, and $w(u,z)=u^{-1}y(\tfrac{1}{2}u,\sqrt {2u} z)$; then (\ref{5.1}) becomes
\begin{equation} \label{5.2}
\frac{d^{2}w}{dz^{2}}+u^{2}(z^{2}+1)w=0.
\end{equation}
As in \cref{sec2} the turning points are at $z=\pm i$ and thus we do not require Airy expansions, since our main interest is the real axis which is free from turning points. The LG expansions will be constructed in the complex plane similarly to those for the parabolic cylinder functions of \cref{sec2}, with $\bar{\xi}$ as given by (\ref{2.2}) and (\ref{2.13}) replaced by $i\bar{\xi}$ in the approximations of this section.

\begin{figure}[htbp]
  \centering
  \includegraphics[width=\textwidth,keepaspectratio]{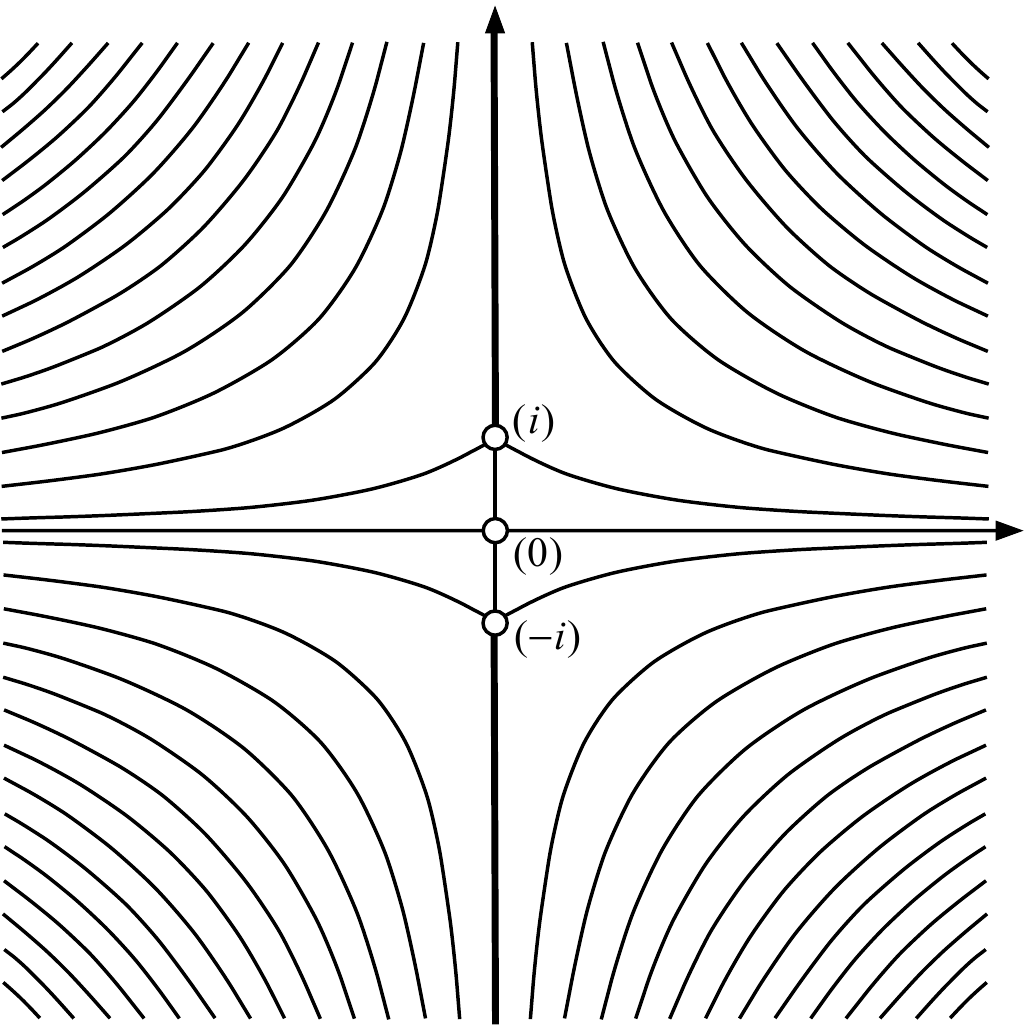}
  \caption{Level curves $\Im(\bar{\xi})=\mathrm{constant}$}
  \label{fig:fig4pcf}
\end{figure}

The level curves determining the regions of validity are thus $\Im(\bar{\xi})=\mathrm{constant}$, and several are depicted in \cref{fig:fig4pcf}. These are orthogonal to those used in \cref{sec2} (\cref{fig:fig1pcf}).

As $z \rightarrow +\infty$ we use (\ref{2.19}) and so similarly to (\ref{4.24}) we arrive at
\begin{multline} \label{5.6}
W_{j}\left(-\tfrac{1}{2}u,\sqrt {2u} z\right)
=\frac{e^{\mp i (\hat{\chi}_{m}(u)
-\frac{1}{8}\pi)}e^{\pi u/8}}
{(2u)^{1/4}\left(1+z^{2}\right)^{1/4}}
\exp \left\{\sum\limits_{s=1}^{m}(-1)^{s}
\frac{\bar{\mathrm{E}}_{2s}(\bar{\beta})}
{u^{2s}}
\right\}
\\
\times \exp \left\{\pm iu\bar{\xi}
\mp i\sum\limits_{s=0}^{m}(-1)^{s}
\frac{\bar{\mathrm{E}}_{2s+1}(\bar{\beta})}
{u^{2s+1}}\right\}
 \left\{1+\bar{\eta}_{2m+2,j}(u,z) \right\},
\end{multline}
where upper signs are for $j=0$ and lower signs for $j=3$. Here
\begin{equation} \label{5.7}
\hat{\chi}_{m}(u)=
\frac{u}{4}\ln\left(\frac{2e}{u}\right)
- \sum\limits_{s=0}^{m}
(-1)^{s}\frac{
\bar{\mathrm{E}}_{2s+1}(1)}
{u^{2s+1}}.
\end{equation}

Error bounds are furnished by (\ref{2.16}) and (\ref{2.17}) with both lower integration limits as $\bar{\beta}=1$, but this value corresponding to $z=e^{\pi i/4}\infty$ for $j=0$ and $z=e^{-\pi i/4}\infty$ for $j=3$. Also as $p$ passes along the paths of integration $\Im (\bar{\xi})$ must be monotonic rather than $\Re(\bar{\xi})$.

The error terms $\bar{\eta}_{2m+2,j}(u,z)$ are $\mathcal{O}(u^{-2m-2})$ uniformly in the right half plane, except in the neighbourhoods of $z=\pm i$. These approximations are also valid in part of left half plane that includes the negative real $z$ axis, but this extension is not required. 

For the derivatives we find in a similar manner to (\ref{2.29}) and (\ref{2.30}) that
\begin{multline} \label{5.8}
W_{j}^{\prime}\left(-\tfrac{1}{2}u,\sqrt {2u} z\right)
=\frac{1}{2} e^{\mp i (\tilde{\chi}_{m}(u)
-\frac{5}{8}\pi)} (2u)^{1/4}
e^{\pi u/8}
\left(1+z^{2}\right)^{1/4}
\\   \times
\exp \left\{
\sum\limits_{s=1}^{m}(-1)^{s}
\frac{\tilde{\mathrm{E}}_{2s}(\bar{\beta})}
{u^{2s}} \right\} \exp \left\{\pm iu\bar{\xi}
\mp i\sum\limits_{s=0}^{m}(-1)^{s}
\frac{\tilde{\mathrm{E}}_{2s+1}(\bar{\beta})}
{u^{2s+1}}\right\}
 \left\{1+\tilde{\eta}_{2m+2,j}(u,z) \right\},
\end{multline}
where
\begin{equation} \label{5.9}
\tilde{\chi}_{m}(u)=
\frac{u}{4}\ln\left(\frac{2e}{u}\right)
-  \sum\limits_{s=0}^{m}
(-1)^{s}\frac{
\tilde{\mathrm{E}}_{2s+1}(1)}
{u^{2s+1}}.
\end{equation}
Error bounds are the same as for the corresponding ones in (\ref{5.6}) except with $\bar{\mathrm{E}}$ replaced by $\tilde{\mathrm{E}}$.

From (\ref{5.8}) and \cite[Eq. 12.2.7]{NIST:DLMF} and using $\tilde{E}_{2s+1}(0)=0$ we see that
\begin{multline} \label{5.10}
\frac{W_{0}^{\prime}
\left(-\tfrac{1}{2}u,0\right)}
{W_{3}^{\prime}
\left(-\tfrac{1}{2}u,0\right)}
= \frac{U^{\prime}\left(-\tfrac{1}{2}iu,
0\right)}
{U^{\prime}\left(\tfrac{1}{2}iu,0\right)}
=-\frac{i2^{iu/2}\Gamma\left(
\tfrac{1}{4}+\tfrac{1}{4}iu\right)}
{\Gamma\left(
\tfrac{1}{4}-\tfrac{1}{4}iu\right)}
\\
=e^{- i (2 \tilde{\chi}_{m}(u)
-\frac{5}{4}\pi)}\left\{
\frac{1+\tilde{\eta}_{2m+2,0}(u,0)}
{1+\tilde{\eta}_{2m+2,3}(u,0)}\right\}.
\end{multline}
Now we get using (\ref{5.10}), similarly to (\ref{4.25}),
\begin{multline} \label{5.11}
e^{-i\left(\frac{1}{2}\phi_{2}(u)+\tilde{\chi}_{m}(u)\right)}
=\left\{\tanh\left(
\tfrac{1}{4}\pi (u+i)\right)\right\}^{1/4}
\left\{\frac{1+\tilde{\eta}_{2m+2,3}(u,0)}
{1+\tilde{\eta}_{2m+2,0}(u,0)}
\right\}^{1/2}
\\ =1+\mathcal{O}(u^{-2m-2}).
\end{multline}
From (\ref{4.27}) this shows $\tilde{\chi}_{m}(u)-\chi_{m}(u)=\mathcal{O}(u^{-2m-2})$ and hence by comparing (\ref{4.22}) and (\ref{5.9}) that the series in both must be equivalent, and hence $\tilde{\chi}_{m}(u)=\chi_{m}(u)$. Thus from (\ref{4.27}) we have
\begin{equation} \label{5.11a}
e^{i\left(\frac{1}{2}\phi_{2}(u)+\tilde{\chi}_{m}(u)\right)}
=e^{i \varepsilon_{m}(u)}
=1+\mathcal{O}(u^{-2m-2}).
\end{equation}

Next, from (\ref{4.3a}) with $a=-u/2$ we have
\begin{multline} \label{5.12}
\bar{k}^{-1/2}(u)W\left(-\tfrac{1}{2}u,\sqrt {2u} x\right)
+i\bar{k}^{1/2}(u)W\left(-\tfrac{1}{2}u,-\sqrt {2u} x\right) \\
=2^{1/2}e^{-\pi u/8}
e^{-i\left(\frac{1}{2}\phi_{2}(u)
-\frac{1}{8}\pi\right)}
W_{0}\left(-\tfrac{1}{2}u,\sqrt {2u}x \right),
\end{multline}
where
\begin{equation} \label{5.13}
\bar{k}(u)=k(-u)=\sqrt{1+e^{-\pi u}} -e^{-\pi u/2}.
\end{equation}
Therefore, from (\ref{5.6}) with $j=0$, (\ref{5.11a}) and (\ref{5.12}), and on expressing
\begin{equation} \label{5.13a}
1+\bar{\eta}_{2m+2,0}(u,x)
=\left\vert 1+\bar{\eta}_{2m+2,0}(u,x) \right\vert
\exp\left(i\arg\left\{1+\bar{\eta}_{2m+2,0}(u,x)
\right\}\right),
\end{equation}
where here $\arg$ is the principal argument, we deduce that
\begin{multline} \label{5.14}
W\left(-\tfrac{1}{2}u,\sqrt {2u} x\right)
=
\left\{\frac{2\bar{k}^{2}(u)}{u\left(1+z^{2}\right)}
\right\}^{1/4}
\exp \left\{\sum\limits_{s=1}^{m}(-1)^{s}
\frac{\bar{\mathrm{E}}_{2s}(\bar{\beta})}
{u^{2s}}\right\}
\\
\times \cos \left\{
u\bar{\xi}+\frac{1}{4} \pi
-\sum\limits_{s=0}^{m}(-1)^{s}
\frac{\bar{\mathrm{E}}_{2s+1}(\bar{\beta})}
{u^{2s+1}}+\tilde{\epsilon}_{2m+2,0}(u,x)
\right\}
 \left\vert 1+\bar{\eta}_{2m+2,0}(u,x) \right\vert,
\end{multline}
and
\begin{multline} \label{5.15}
W\left(-\tfrac{1}{2}u,-\sqrt {2u} x\right)
=\left\{
\frac{2}{u \bar{k}^{2}(u)\left(1+z^{2}\right)}
\right\}^{1/4}\exp \left\{
\sum\limits_{s=1}^{m}(-1)^{s}
\frac{\bar{\mathrm{E}}_{2s}(\bar{\beta})}
{u^{2s}}\right\}
\\
\times \sin \left\{
u\bar{\xi} +\frac{1}{4} \pi
-\sum\limits_{s=0}^{m}(-1)^{s}
\frac{\bar{\mathrm{E}}_{2s+1}(\bar{\beta})}
{u^{2s+1}}+\tilde{\epsilon}_{2m+2,0}(u,x)
\right\}
 \left\vert 1+\bar{\eta}_{2m+2,0}(u,x) \right\vert,
\end{multline}
where
\begin{equation} \label{5.16}
\tilde{\epsilon}_{2m+2,0}(u,x)
=\arg\left\{1+\bar{\eta}_{2m+2,0}(u,x)
\right\}-\varepsilon_{m}(u).
\end{equation}
Recall that $\varepsilon_{m}(u)$ is explicitly given by (\ref{4.27}) and is $\mathcal{O}(u^{-2m-2})$ as $u \rightarrow \infty$. In addition $\bar{\eta}_{2m+2,0}(u,x)=\mathcal{O}(u^{-2m-2})$, uniformly for $0\leq x < \infty$, and hence the same is true for $\tilde{\epsilon}_{2m+2,0}(u,x)$. We emphasize the double asymptotic behaviour that $\bar{\eta}_{2m+2,0}(u,x) \rightarrow 0$ as $x \rightarrow \infty$.

\subsection{Inhomogeneous equation}

Solutions of the inhomogeneous equation
\begin{equation} \label{5.17}
\frac{d^{2}y}{dz^{2}}+\left(\frac{1}{4}z^{2}+a\right)y
=z^{R},
\end{equation}
are given by $W_{R}^{(j,k)}(-a,z)$. With the same change of variable and parameter as in previous sections this is transformed to
\begin{equation} \label{5.20}
\frac{d^{2}w}{dz^{2}}+u^{2}(z^{2}+1)w=z^R,
\end{equation}
and for  $j,k\in \{0,1,2,3\}$ with $j<k$ asymptotic solutions are furnished by
\begin{equation} \label{5.18}
\tilde{w}_{R}^{(j,k)}(u,z)=\frac{1}{u^{2}}\sum\limits_{s=0}^{n-1} (-1)^{s+1}\frac{\bar{G}_{s,R}(z)}{u^{2s}} +\tilde{\varepsilon}_{n,R}^{(j,k)}(u,z),
\end{equation}
where $\bar{G}_{s,R}(z)$ are given by (\ref{2.41}) and (\ref{2.42}). The error terms satisfy (\ref{2.43}) and (\ref{2.44}) where the paths of integration are as described for those bounds, except that they run from $z=e^{\pi i/4}i^{j}\infty$ to $z=e^{\pi i/4}i^{k}\infty$, and the monotonicity condition is for $\Im(\bar{\xi})$.

From (\ref{4.35a}) a particular solution is given by
\begin{multline} \label{5.21a}
 W_{R}^{(0,3)}(-a,z)=ie^{-\pi a/2}
\left[W_{3}(-a,z) \int_{e^{\pi i/4}\infty}^z \right.
 t^{R}W_{0}(-a,t)dt \\
\left. -W_{0}(-a,z) \int_{e^{-\pi i/4}\infty}^z 
 t^{R}W_{3}(-a,t)dt \right].
\end{multline}
Again this is the solution of most interest, since it is not exponentially large in the right half plane. In fact, as we shall show, unlike $ W_{R}^{(0,3)}(a,z)$ it is also not exponentially large in part of the left half plane including the whole negative real axis. Matching this solution with corresponding asymptotic solution that is also not exponentially large in the right half plane yields
\begin{equation} \label{5.21}
W_{R}^{(0,3)}\left(-\tfrac{1}{2}u,\sqrt{2u} z\right)
=(2u)^{\frac{1}{2}R+1}\tilde{w}_{R}^{(-1,1)}(u,z).
\end{equation}

\begin{figure}[htbp]
  \centering
  \includegraphics[width=\textwidth,keepaspectratio]{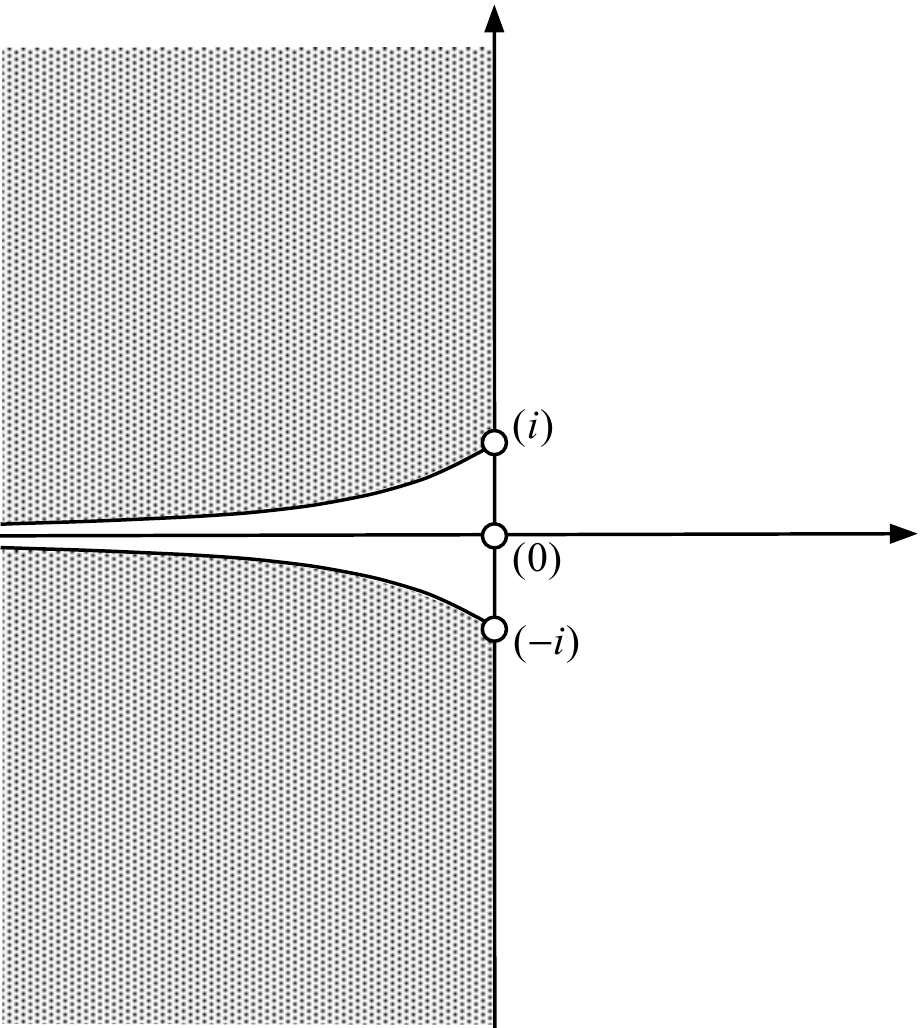}
  \caption{Domain $\breve{Z}^{(0,3)}$}
  \label{fig:fig4apcf}
\end{figure}

Note that this solution is real for real $z$. See \cref{fig:fig4apcf} for region of validity $\breve{Z}^{(0,3)}$, all boundaries of which must be excluded: these are the parts of the imaginary axis from $z=\pm i$ to $z=\pm i\infty$, and the level curves in the left half plane emanating from $z=\pm i$. Note the region includes the whole real $z$ axis, with the error term being $\mathcal{O}(u^{-2n-2})$ as $u \rightarrow \infty$ and also is vanishing as $z \rightarrow \infty$. However, it does not vanish as $z \rightarrow -\infty$.

To get a sharper approximation for $\Re(z) \leq 0$ we use a suitable connection formula, derived as follows. Firstly, from (\ref{4.40})
\begin{multline}  \label{5.22}
W_{R}^{(0,2)}\left(-\tfrac{1}{2}u,\sqrt{2u} z\right)
=W_{R}^{(2,3)}\left(-\tfrac{1}{2}u,\sqrt{2u} z\right)
\\ +ie^{-3 R\pi i/4}\Lambda_{R}\left(-\tfrac{1}{2}iu\right)
W_{2}\left(-\tfrac{1}{2}u,\sqrt{2u} z\right).
\end{multline}

Next we can assert there is a constant $\alpha_{R}(u)$ such that
\begin{equation} \label{5.23}
W_{R}^{(0,3)}\left(-\tfrac{1}{2}u,\sqrt{2u} z\right)
=W_{R}^{(0,2)}\left(-\tfrac{1}{2}u,\sqrt{2u} z\right)
+\alpha_{R}(u) W_{0}\left(-\tfrac{1}{2}u,\sqrt{2u} z\right),
\end{equation}
since the difference of the two particular solutions of (\ref{5.20}) appearing in this equation is a solution of the homogeneous equation (\ref{5.2}) which, like both of them, is not exponentially large in the first quadrant; this difference must then be a multiple of $W_{0}\left(-\tfrac{1}{2}u,\sqrt{2u} z\right)$.

We can find $\alpha_{R}(u)$ by letting $z \rightarrow \infty e^{-\pi i/4}$. In this limit we have from (\ref{5.23})
\begin{equation} \label{5.24}
W_{R}^{(0,2)}\left(-\tfrac{1}{2}u,\sqrt{2u} z\right)
\sim -\alpha_{R}(u) W_{0}\left(-\tfrac{1}{2}u,\sqrt{2u} z\right),
\end{equation}
since in this case $W_{R}^{(0,3)}\left(-\tfrac{1}{2}u,\sqrt{2u} z\right)$ is $\mathcal{O}(z^{R-2})$ and the other two are exponentially large. On the other hand, in the same limit we have from (\ref{5.22})
\begin{equation} \label{5.25}
W_{R}^{(0,2)}\left(-\tfrac{1}{2}u,\sqrt{2u} z\right)
\sim ie^{-3 R\pi i/4}
\Lambda_{R}\left(-\tfrac{1}{2}iu\right)
W_{2}\left(-\tfrac{1}{2}u,\sqrt{2u} z\right).
\end{equation}
Therefore from (\ref{4.3}), (\ref{5.23}), (\ref{5.24}) and (\ref{5.25})
\begin{multline} \label{5.26}
\alpha_{R}(u)= -ie^{-3 R\pi i/4}
 \Lambda_{R} \left(-\tfrac{1}{2}iu\right)
 \lim_{z \rightarrow \infty \exp(-\pi i/4)}
 \frac{W_{2}\left(-\tfrac{1}{2}u,\sqrt{2u} z\right)}{W_{0}\left(-\tfrac{1}{2}u,\sqrt{2u} z\right)}
 \\=-e^{-3 R\pi i/4}e^{-\pi u/2}
 \Lambda_{R}\left(-\tfrac{1}{2}iu\right),
\end{multline}
and hence from (\ref{2.55})
\begin{multline} \label{5.27}
 \alpha_{R}(u)=\sqrt {\pi }e^{-\frac{1}{4}\pi u+\frac{1}{4}(1-R)\pi i}
 2^{\frac{3}{2}R+\frac{1}{4}+\frac{1}{4}iu}
\Gamma \left(\tfrac{1}{2}-\tfrac{1}{2}iu \right) \\
\times
\mathbf{F}\left(\tfrac{1}{2}-\tfrac{1}{2}R,
-\tfrac{1}{2}R;\tfrac{3}{4}
-\tfrac{1}{2}R-\tfrac{1}{4}iu;\tfrac{1}{2} \right).
\end{multline}

Next, from (\ref{5.23}) we obtain by replacing $z$ by $-z$
\begin{equation} \label{5.28}
W_{R}^{(0,3)}\left(-\tfrac{1}{2}u,-\sqrt{2u} z\right)
=(-1)^{R}W_{R}^{(0,2)}\left(-\tfrac{1}{2}u,\sqrt{2u} z\right) +\alpha_{R}(u) W_{0}\left(-\tfrac{1}{2}u,-\sqrt{2u} z\right),
\end{equation}  since by uniqueness $W_{R}^{(0,2)}(-a,-z)=(-1)^{R}W_{R}^{(0,2)}(-a,z)$, both being solutions of (\ref{5.20}) that are not exponentially large in the first and third quadrants. So from (\ref{5.23}) and (\ref{5.28})
\begin{multline} \label{5.29}
W_{R}^{(0,3)}\left(-\tfrac{1}{2}u,\sqrt{2u} z\right)
=(-1)^{R}W_{R}^{(0,3)}\left(-\tfrac{1}{2}u,-\sqrt{2u} z\right)  \\
+(-1)^{R+1}\alpha_{R}(u) W_{0}\left(-\tfrac{1}{2}u,-\sqrt{2u} z\right) 
+\alpha_{R}(u) W_{0}\left(-\tfrac{1}{2}u,\sqrt{2u} z\right).
\end{multline}

We then finally use
\begin{equation} \label{5.30}
W_{0}\left(-\tfrac{1}{2}u,\sqrt{2u} z\right)=
c_{0}(u)W_{0}\left(-\tfrac{1}{2}u,-\sqrt{2u} z\right)
+c_{3}(u)W_{3}\left(-\tfrac{1}{2}u,-\sqrt{2u} z\right),
\end{equation}
where from (\ref{4.2}) and \cite[Eq. 12.2.17]{NIST:DLMF}
\begin{equation} \label{5.32}
c_{0}(u)=-ie^{-\pi u/2},
\end{equation}
and
\begin{equation} \label{5.31}
c_{3}(u)=
\frac{\sqrt{2 \pi}e^{\pi i/4}e^{-\pi u/4}}
{\Gamma  \left( \frac{1}{2}-\frac{1}{2}iu \right)}.
\end{equation}
Consequently, from (\ref{5.29}), (\ref{5.30}) and (\ref{5.32}) we get our desired connection formula
\begin{multline} \label{5.33}
W_{R}^{(0,3)}\left(-\tfrac{1}{2}u,\sqrt {2u} z\right)
=(-1)^{R}W_{R}^{(0,3)}\left(-\tfrac{1}{2}u,-\sqrt {2u} z\right)\\
+\left\{(-1)^{R+1}-ie^{-\pi u/2}
\right\}\alpha_{R}(u)
W_{0}\left(-\tfrac{1}{2}u,-\sqrt {2u} z\right)
+\alpha_{R}(u) c_{3}(u)W_{3}\left(-\tfrac{1}{2}u,-\sqrt {2u} z\right),
\end{multline}
where $\alpha_{R}(u)$ is given by (\ref{5.27}) and $c_{3}(u)$ by (\ref{5.31}). 

We can go a step further, noting that when $z$ is real so too is $W_{R}^{(0,3)}(a,z)$, and moreover $W_{0}(a,z)$ and $W_{3}(a,z)$ are complex conjugates. Hence the coefficients of these two functions in (\ref{5.33}) must be also complex conjugates. This allows us to recast this connection formula into the form
\begin{multline} \label{5.34}
W_{R}^{(0,3)}\left(-\tfrac{1}{2}u,\sqrt {2u} z\right)
=(-1)^{R}W_{R}^{(0,3)}\left(-\tfrac{1}{2}u,-\sqrt {2u} z\right)  \\
+\left\{(-1)^{R+1}-ie^{-\pi u/2}
\right\}\alpha_{R}(u)
W_{0}\left(-\tfrac{1}{2}u,-\sqrt {2u} z\right)  \\
+\left\{(-1)^{R+1}+ie^{-\pi u/2}
\right\}\overline{\alpha_{R}(u)}W_{3}\left(-\tfrac{1}{2}u,-\sqrt {2u} z\right),
\end{multline}
where $\overline{\alpha_{R}(u)}$ is given by (\ref{5.27}) with $i$ replaced by $-i$.

For $\Re(z) \leq 0$ we can use this, along with (\ref{5.6}), (\ref{5.18}) and (\ref{5.21}), to obtain an asymptotic expansion which is uniformly valid in this left half plane except near $z=\pm iy$, $1 \leq y < \infty$. We also note that all error terms vanish as $\Re(z) \rightarrow -\infty$.

We observe that the last two terms on the RHS of (\ref{5.34}) are exponentially small in the unshaded region of the left half plane in \cref{fig:fig4apcf}, and hence the first term, which has the asymptotic expansion (\ref{5.18}), dominates here (as expected). As $z$ enters the shaded region in the second quadrant the second term on the RHS term becomes exponentially large and it dominates, and as $z$ enters the shaded region in the the third quadrant the third term on the RHS term becomes exponentially large and it dominates. This is an example of the Stokes phenomenon where a function's asymptotic behaviour changes from subdominant to dominant, or vice versa, as its argument crosses a Stokes line.

\section*{Acknowledgments}
Financial support from Ministerio de Ciencia e Innovaci\'on, Spain, project PGC2018-098279-B-I00 (MCIU/AEI/FEDER, UE) is acknowledged. 

\bibliographystyle{siamplain}
\bibliography{biblio}

\begin{thebibliography}{10}

\bibitem{Zaytoon:2016:WID}
{\sc M.~Abu~Zaytoon, T.~Alderson, and M.~Hamdan}, {\em Weber's inhomogeneous
  differential equation with initial and boundary conditions}, Int. J. Open
  Problems Compt. Math, 9 (2016), pp.~1--11,
  \url{https://doi.org/10.12816/0033917}.

\bibitem{Buchholz:1969:CHF}
{\sc H.~Buchholz}, {\em The confluent hypergeometric function}, Springer,
  Berlin, Heidelberg, 1969, \url{https://doi.org/10.1007/978-3-642-88396-5}.

\bibitem{NIST:DLMF}
{\em {\it NIST Digital Library of Mathematical Functions}}.
\newblock Release 1.0.25 of 2019-12-15, \url{http://dlmf.nist.gov/}.
\newblock F.~W.~J. Olver, A.~B. {Olde Daalhuis}, D.~W. Lozier, B.~I. Schneider,
  R.~F. Boisvert, C.~W. Clark, B.~R. Miller, B.~V. Saunders, H.~S. Cohl, and
  M.~A. McClain, eds.

\bibitem{Dunster:2020:ASI}
{\sc T.~M. Dunster}, {\em Asymptotic solutions of inhomogeneous differential
  equations having a turning point.}, Stud. Appl. Math., 145 (2020),
  pp.~500--536, \url{https://doi.org/10.1111/sapm.12326}.

\bibitem{Dunster:2020:LGE}
{\sc T.~M. Dunster}, {\em Liouville-{G}reen expansions of exponential form,
  with an application to modified {B}essel functions}, Proc. Roy. Soc.
  Edinburgh Sec. A, 150 (2020), pp.~1289--1311,
  \url{https://doi.org/10.1017/prm.2018.117}.

\bibitem{Dunster:2017:COA}
{\sc T.~M. Dunster, A.~Gil, and J.~Segura}, {\em Computation of asymptotic
  expansions of turning point problems via {C}auchy's integral formula: Bessel
  functions.}, Constr. Approx., 46 (2017), pp.~645--675,
  \url{https://doi.org/10.1007/s00365-017-9372-8}.

\bibitem{Dunster:2018:USE}
{\sc T.~M. Dunster, A.~Gil, and J.~Segura}, {\em Uniform asymptotic expansions
  for {L}aguerre polynomials and related confluent hypergeometric functions},
  Advances in Computational Mathematics, 44 (2018), pp.~1441--1474,
  \url{https://doi.org/10.1007/s10444-018-9589-5}.

\bibitem{Dunster:2020:SEB}
{\sc T.~M. Dunster, A.~Gil, and J.~Segura}, {\em Simplified error bounds for
  turning point expansions}, Anal. Appl.,  (2020),
  \url{https://doi.org/10.1142/S0219530520500104}.

\bibitem{Gil:2004:IRC}
{\sc A.~Gil, J.~Segura, and N.~M. Temme}, {\em Integral representations for
  computing real parabolic cylinder functions}, Numer. Math., 98 (2004),
  pp.~105--134, \url{https://doi.org/10.1007/s00211-004-0517-x}.

\bibitem{Gil:2006:RPC}
{\sc A.~Gil, J.~Segura, and N.~M. Temme}, {\em Algorithm 850: Real parabolic
  cylinder functions {$U(a,x)$}, {$V(a,x)$}}, ACM Trans. Math. Software, 32
  (2006), pp.~102--112, \url{https://doi.org/10.1145/1132973.1132978}.

\bibitem{Gil:2006:CRPC}
{\sc A.~Gil, J.~Segura, and N.~M. Temme}, {\em Computing the real parabolic
  cylinder functions {$U(a,x)$}, {$V(a,x)$}}, ACM Trans. Math. Software, 32
  (2006), pp.~70--101, \url{https://doi.org/10.1145/1132973.1132977}.

\bibitem{Gil:2011:A914}
{\sc A.~Gil, J.~Segura, and N.~M. Temme}, {\em Algorithm 914: parabolic
  cylinder function {$W(a,x)$} and its derivative}, ACM Trans. Math. Software,
  38 (2011), pp.~Art. 6, 5, \url{https://doi.org/10.1145/2049662.2049668}.

\bibitem{Gil:2011:FAC}
{\sc A.~Gil, J.~Segura, and N.~M. Temme}, {\em Fast and accurate computation of
  the {W}eber parabolic cylinder function {$W(a,x)$}}, IMA J. Numer. Anal., 31
  (2011), pp.~1194--1216, \url{https://doi.org/10.1093/imanum/drq012}.

\bibitem{Jones:2006:PCF}
{\sc D.~S. Jones}, {\em Parabolic cylinder functions of large order}, J.
  Comput. Appl. Math., 190 (2006), pp.~453--469,
  \url{https://doi.org/10.1016/j.cam.2005.04.009}.

\bibitem{Nield:2009:TET}
{\sc D.~A. Nield and A.~V. Kuznetsov}, {\em The effect of a transition layer
  between a fluid and a porous medium: shear flow in a channel}, Transp. Porous
  Media, 78 (2009), pp.~477--487.

\bibitem{Olver:1959:USE}
{\sc F.~W.~J. Olver}, {\em Uniform asymptotic expansions for {W}eber parabolic
  cylinder functions of large orders}, J. Res. Nat. Bur. Stand., 63B (1959),
  pp.~131--169.

\bibitem{Olver:1975:SOL}
{\sc F.~W.~J. Olver}, {\em Second-order linear differential equations with two
  turning points}, Philos. Trans. R. Soc. A, 278 (1975), pp.~137--174,
  \url{https://doi.org/10.1098/rsta.1975.0023}.

\bibitem{Olver:1997:ASF}
{\sc F.~W.~J. Olver}, {\em Asymptotics and special functions}, AKP Classics, A
  K Peters Ltd., Wellesley, MA, 1997.
\newblock Reprint of the 1974 original [Academic Press, New York].

\bibitem{Prudnikov:1986:IAS}
{\sc A.~P. Prudnikov, I.~U.~A. Brychkov, and O.~I. Marichev}, {\em Integrals
  and Series: Special functions}, Integrals and Series, Gordon and Breach
  Science Publishers, 1986.

\bibitem{Temme:2000:NAA}
{\sc N.~M. Temme}, {\em Numerical and asymptotic aspects of parabolic cylinder
  functions}, J. Comput. Appl. Math., 121 (2000), pp.~221--246,
  \url{https://doi.org/10.1016/S0377-0427(00)00347-2}.

\bibitem{Temme:2015:AMF}
{\sc N.~M. Temme}, {\em Asymptotic methods for integrals}, vol.~6 of Series in
  Analysis, World Scientific Publishing Co. Pte. Ltd., Hackensack, NJ, 2015.

\end{thebibliography}

\end{document}